\documentclass[10pt]{amsart}
\usepackage{amsmath,amsfonts,amssymb,color,a4,epsfig,graphics}
\usepackage[english]{babel}
\usepackage{enumerate, mathrsfs}
\usepackage{xcolor}
\usepackage{tikz}
\def\build#1_#2^#3{\mathrel{\mathop{\kern 0pt#1}\limits_{#2}^{#3}}}
\usetikzlibrary{arrows}

\newcommand{\R}{{\mathbb{R}}}
\newcommand{\C}{{\mathbb{C}}}
\newcommand{\Z}{{\mathbb{Z}}}
\newcommand{\N}{{\mathbb{N}}}

\newcommand{\Ac}{\mathcal{A}}
\newcommand{\Bc}{\mathcal{B}}
\newcommand{\Cc}{\mathcal{C}}
\newcommand{\Dc}{\mathcal{D}}
\newcommand{\Ec}{\mathcal{E}}
\newcommand{\Fc}{\mathcal{F}}
\newcommand{\Gc}{\mathcal{G}}
\newcommand{\Hc}{\mathcal{H}}
\newcommand{\Ic}{\mathcal{I}}
\newcommand{\Kc}{\mathcal{K}}
\newcommand{\Nc}{\mathcal{N}}

\newcommand{\Vc}{\mathcal{V}}
\newcommand{\Dr}{\mathscr{D}}
\newcommand{\un}{{\rm \bf {1}}}
\newcommand{\supp}{{\rm supp }\,}

\numberwithin{equation}{section}

\def\id{\mathop{\rm id}\nolimits}

\def\rmi{{\rm i}}
\def \lint{[\![}
\def \rint{]\!]}

\newtheorem{theorem}{Theorem}[section]
\newtheorem{proposition}[theorem]{Proposition}
\newtheorem{lemma}[theorem]{Lemma}
\newtheorem{remark}[theorem]{Remark}
\newtheorem{definition}[theorem]{Definition}
\newtheorem{corollary}[theorem]{Corollary}

\begin{document}

\title[Spectral analysis of the discrete Laplacian]{Spectral analysis of the Laplacian acting
  on discrete cusps and funnels}
\author{Nassim Athmouni}
\address{Universit\'e de Gafsa, Campus Universitaire 2112, Tunisie}
\email{\tt athmouninassim@yahoo.fr}
 \author{Marwa Ennaceur}
 \address{Universit\'e de Sfax. Route de la Soukra km 4 - B.P. n° 802 - 3038 Sfax}
\email{\tt ennaceur.marwa27@gmail.com}
\author{Sylvain Gol\'enia}
\address{Univ. Bordeaux,
Bordeaux INP, CNRS, IMB, UMR 5251, F-33400 Talence, France}
\email{\tt sylvain.golenia@math.u-bordeaux.fr}

\begin{abstract}We study perturbations of the discrete Laplacian associated to discrete analogs of cusps and funnels. We perturb the metric and the potential in a long-range way. We establish a propagation estimate and a Limiting Absorption Principle away from the possible embedded eigenvalues. The approach is based on a positive commutator technique.
\end{abstract}

\maketitle
\tableofcontents
\section{Introduction}
The spectral theory of discrete Laplacians on graphs has drawn a lot
of attention for decades as they are discrete analogs of manifolds. 
We are especially interested in the nature of the essential spectrum. Without trying to be exhaustive, using positive commutator techniques, \cite{S,BoSa} treat the case of $\Z^d$, \cite{AF, GG} study the case of binary trees, \cite{MRT} investigate some general graphs, and \cite{PR} focused on a periodic setting. Some other techniques have been used successfully, e.g., \cite{HN} with some geometric approach and \cite{BrKe}.

In the context of some manifolds of finite volume, \cite{AbTr,GoMo} prove
that the essential spectrum of the (continuous) Laplacian becomes empty under the presence of a magnetic field with compact support. Besides, they establish some Weyl asymptotic. Analogously, for some discrete cusps, \cite{GT} classify magnetic potentials that lead to the absence of the essential spectrum and compute a kind of Weyl asymptotic for the magnetic discrete Laplacian. Back to \cite{GoMo}, one also obtains a refined analysis of the spectral measure (propagation estimate, limiting absorption principle) for long-range perturbation of the metric when the essential spectrum occurs relying on a positive commutator technique. We refer to \cite{GoMo} for further comments and references therein. This part of the analysis was not carried out in \cite{GT}. This is the main aim of this article.


To start off, we recall some standard definitions of graph theory. A (non-oriented) \emph{graph} is
a triple
$\Gc:=(\mathcal{E},\mathcal{V}, m)$, where $\mathcal{V}$ is a finite or countable
set (the \emph{vertices}), $\mathcal{E}:\mathcal{V} \times
\mathcal{V}\rightarrow
\mathbb{R}_{+}$ 
is symmetric, and $m:\Vc\to (0,
\infty)$ is a weight. 
We
say that $\Gc$ is \emph{simple} if $m=1$ and $\Ec:\Vc\times \Vc\to
\{0,1\}$.

Given  $x,y\in \mathcal{V}$, we say that $(x,y)$ is an \emph{edge} and that $x$ and $y$ are
\emph{neighbors} if $\mathcal{E}(x,y)>0$. Note that in this case,  since $\mathcal{E}$ is symmetric, $(y,x)$ is also an edge and $y$ and $x$ are neighbors.  We denote this
relationship by $x\sim y$ and the set of neighbors of $x$ by
$\Nc_\Gc(x)$. The space of complex-valued functions acting on the set of vertices $\Vc$ is
denoted by $C(\Vc):=\lbrace f:\Vc \to {\C}\rbrace$. Moreover,  $C_c(\Vc)$ is the subspace of $C(\Vc)$  of functions with finite  support. We consider the Hilbert space
\[\ell^2(\Vc, m):=\left\{ f\in C(\Vc), \quad \sum_{x\in\Vc} m(x)|f(x)|^2 <\infty \right\}\]
endowed with the scalar product, $\langle f,g\rangle:= \sum_{x\in \Vc}m(x) \overline{f(x)}g(x)$. We define the Laplacian operator
\begin{align}\label{e:DeltaG}
\Delta_{\Gc} f(x):=\frac{1}{m(x)}\sum_{y\in\Vc}\Ec(x,y)(f(x)-f(y)),
\end{align}
for all $f\in\Cc_c(\Vc)$. $\Delta_{\Gc}$ is a positive operator since we have $\langle f, \Delta_{\Gc} f\rangle_{\ell^2(\Vc,m)}= Q_{\Gc} (f)$, with \[ Q_{\Gc} (f):= \frac{1}{2}\sum_{x,y \in \Vc} \Ec(x,y) \left| f(x) -f(y)\right|^2,\]
 for all $f\in \Cc_c(\Vc)$. To simplify, we denote its Friedrichs'extension with the same symbol. We define the \emph{degree} of $x\in \Vc$ by
\[\deg_\Gc(x):=\frac{1}{m(x)}\sum_{y\in \Vc} \Ec(x,y)\]

We present a simple version of our model:  We consider $\Gc_1:=(\Ec_1, \Vc_1, m_1)$, where $\Vc_1:=\Z$, $m_1(n):=e^{-n}$,
and $\Ec(n,n+1):=e^{-(2n+1)/2}$,
for all $n\in\N$ and $\Gc_2:=(\Ec_2,\Vc_2,m_2)$ a connected finite graph such that $|\Vc_2|=p$, $p\geq3$, where $|\Vc_2| $ is the cardinal of the set $\Vc_2$ with $m_2$ constant. Let $\Gc:=(\Ec, \Vc, m)$ be the \emph{twisted cartesian product} $\Gc_1\times_{\tau} \Gc_2$ given by
\begin{align*}\left\{\begin{array}{rl}
m(x,y):=& m_1(x)\times m_2(y),
\\
\Ec\left((x,y),(x',y')\right):=
&\Ec_1(x,x')\times \delta_{y, y'} +\delta_{x,x' } \times \Ec_2(y,y'),
\end{array}\right.\end{align*} for all $x,x'\in\Vc_1$ and $y,y'\in \Vc_2$,
\begin{figure}\label{figure}
\begin{tikzpicture}[scale=0.55]
\def\test{0.75}
\def\xmax{10}
  \foreach \x in {1,2,...,\xmax}
{\fill[color=black]({(\x)/\test}, 0)circle(.3mm);
\fill[color=black]({(\x)/\test+.25/\test/\x}, 4*0.75/\x)circle(.3mm);
\fill[color=black]({(\x)/\test-1/\test/\x}, 4*1/\x)circle(.3mm);
\draw(\x/\test, 0)--(\x/\test+.25/\test/\x, 4*0.75/\x);
\draw(\x/\test, 0)--(\x/\test-1/\test/\x, 4*1/\x);
\draw(\x/\test+.25/\test/\x, 4*0.75/\x)--(\x/\test-1/\test/\x, 4*1/\x);
\draw(\x/\test, 0)--({(\x+1)/\test}, 0);
\draw(\x/\test-1/\test/\x, 4*1/\x)--({(\x+1)/\test-1/\test/(\x+1)}, {4*1/(\x+1)});
\draw(\x/\test+.25/\test/\x, 4*0.75/\x)--({(\x+1)/\test+.25/\test/(\x+1)}, {4*.75/(\x+1)});
}
\fill[color=black]({(\xmax+1)/\test}, 0)circle(0.3mm);
\fill[color=black]({(\xmax+1)/\test+.25/\test/(\xmax+1)}, {4*0.75/(\xmax+1)})circle(.3mm);
\fill[color=black]({(\xmax+1)/\test-1/\test/(\xmax+1)}, {4*1/(\xmax+1)})circle(.3mm);
\draw({(\xmax+1)/\test}, 0)--({(\xmax+1)/\test+.25/\test/(\xmax+1)}, {4*0.75/(\xmax+1)});
\draw({(\xmax+1)/\test}, 0)--({(\xmax+1)/\test-1/\test/(\xmax+1)}, {4*1/(\xmax+1)});
\draw({(\xmax+1)/\test+.25/\test/(\xmax+1)}, {4*0.75/(\xmax+1)})--({(\xmax+1)/\test-1/\test/(\xmax+1)}, {4*1/(\xmax+1)});
\path({(\xmax+1)/\test+1}, 0) node {$\cdots$};
\path(3, -0.5) node {\emph{funnel side }};
\path(13, -0.5) node {\emph{cusp side }};
\end{tikzpicture}
\caption{Representation of a discrete cusp and funnel side}
\end{figure}
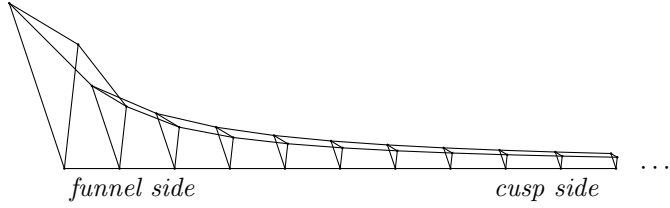
If $n>0$, this is a cups side and if $n<0$, this is a funnel side. We refer to Section \ref{s:setting} for more details.

The (twisted cartesian) Laplacian $\Delta_\Gc$ is essentially self-adjoint on $\Cc_c(\Vc)$, see Proposition \ref{p:cuspessaa}. Moreover, it has no singularly continuous spectrum and
 \[\sigma_{\rm ac}(\Delta_{\Gc})= \left[\frac{\alpha}{m_2},\frac{\beta}{m_2}\right],\]
 with \begin{align}\label{e:c1c2}
\alpha := e^{1/2} + e^{-1/2} -2 \quad \mbox{ and } \quad \beta := e^{1/2} + e^{-1/2} +2.
\end{align}
We turn into perturbation theory. First, we perturb the weights, we consider
$\Gc':=(\Ec',\Vc,m')$, where
 \begin{align*}
m'(x):=(1+\mu(x))m(x) \text{ and } \Ec'(x,y):=(1+\varepsilon(x,y))\Ec(x,y),
 \end{align*}
 \begin{align*}
  (H_0) \hspace*{1.8 cm }\left\{\begin{array}{ll}
 \displaystyle \max_{x_2\in \Vc_2}|V((x_1, x_2))|\to 0, & \text{ if }  |x_1|\rightarrow\infty,
\\
 \displaystyle \max_{x_2\in \Vc_2}|\mu((x_1, x_2))|\rightarrow0,& \text{ if } |x_1|\rightarrow\infty,
 \\
\displaystyle \max_{x_2\in \Vc_2, y\sim (x_1, x_2) } |\varepsilon((x_1, x_2), y)|\rightarrow 0, &\text { if } |x_1|\rightarrow\infty.
 \end{array}\right.
\end{align*}
This ensures that $\Delta_{\Gc'}+V(\cdot)$ is also essentially self-adjoint on $\Cc_c(\Vc)$. Here $V(\cdot)$ denotes the operator of multiplication by $V$. Moreover, $(H_0)$ guarantees the stability of the essential spectrum, see Proposition \ref{com}. Namely,
 \[\sigma_{\rm ess} (\Delta_{\Gc'})= \left[\frac{\alpha}{m_2},\frac{\beta}{m_2}\right].\]
 In order to obtain the absence of singularly continuous spectrum for $\Delta_{\Gc'}$,
 we require some additional decay. Let $\epsilon>0$ and ask:
 \begin{align*}
 (H_1)\quad& \sup_{n\in\Z, y\in \Vc_2} \langle n\rangle^{1+\epsilon}  \Big|V(n-1,y)-V(n,y)\Big| <\infty,
 \\
 (H_2)\quad& \sup_{n\in\Z, y\in \Vc_2} \langle n\rangle^{1+\epsilon}\left|\mu(n-1, y)-\mu(n, y)\right|<\infty
 \\
 (H_3)\quad& \sup_{n\in \Z, k\in \Vc_2} \langle n\rangle^{1+\epsilon}\left|\varepsilon((n,k), (n+1,k))-\varepsilon((n-1,k),(n,k))\right|<\infty,
 \end{align*}
where $\langle \cdot \rangle := \sqrt{1+|\cdot|^2}$.


Our main result is the following:
\begin{theorem}\label{t:LAP} Let $H:= \Delta_{\Gc'} + V(\cdot)$ as above. Suppose that $(H_0)$ holds true.
Then, we have the following assertions:
\begin{enumerate}
\item[(1)] $\sigma_{ess}(H)=\sigma_{ess}(\Delta_{\Gc})$.
\end{enumerate}
Assume furthermore that $(H_1), (H_2)$, and $(H_3)$ hold true. Set $\kappa(H):=\sigma_{p}(H)\cup\{\alpha/m_2,\beta/m_2\}$
with $\alpha,\beta$ are given in \eqref{e:c1c2}
and where $\sigma_p$ denotes the pure point spectrum. Take $s>1/2$ and $[a,b]\subset \R\setminus \kappa(H)$. We obtain:
\begin{enumerate}
\item[(2)] The eigenvalues of $H$ distinct from $\alpha/m_2$ and $\beta/m_2$ are of finite multiplicity and can accumulate
only toward $\alpha/m_2$ and $\beta/m_2$.
\item[(3)] The singular continuous spectrum of $H$  is empty.
\item[(4)] The following limit exists and finite:
\[\lim_{\rho\to0}  \sup_{\lambda \in [a,b]}\|\langle \Lambda\rangle^{-s}(H-\lambda-\rmi\rho)^{-1}\langle \Lambda\rangle^{-s}\|<\infty,\]
\item[(5)] There exists $c>0$ such that for all $f\in\ell^2(\Vc, m')$,  we have:
\[\int_{\R}\|\langle \Lambda\rangle^{-s}e^{-\rmi tH}E_{[a,b]}(H)f\|^2_{\ell^2(\Vc, m')}dt\leq c\|f\|^2_{\ell^2(\Vc, m')}.\]
\end{enumerate}
\end{theorem}
Our approach is based on a positive commutator technique, namely we establish a Mourre estimate. The proof of this theorem is given in Subsection \ref{P.result}. We refer to Section \ref{s:mourre} for historical references and for an introduction on the subject.

We now describe the structure of the paper. In Section \ref{s:mourre}, we present the Mourre's theory. The next section is devoted to study the free model. In Subsection \ref{s:setting}, we present the context and introduce the notion of cusp and funnel. In Subsection \ref{Mourre N}, we start with the Mourre estimate on $\N$. In Subsections \ref{section3}, and \ref{section4}, we prove the  Mourre estimate for the unperturbed Laplacian that acts on a funnel and on a cusp, respectively. Then, in Subsection \ref{section6}, we conclude the Mourre estimate for the whole graph. In Section \ref{section7}, we perturb the metrics and add a potential. The proofs are more involved than in Section \ref{s:setting} as we rely on the optimal class  $\Cc^{1,1}(\Ac)$ of the Mourre theory. This yields the main result.

\noindent{{\bf Notation:} We denote by $\N$ the set of non-negative integers. In particular $0\in \N$. Set $\lint a, b\rint:= [a,b]\cap \Z$.  We denote by $\un_X$ the indicator of the set $X$. }

\noindent{\bf Acknowledgements: } We would like to anonymous referees for their comments on the script.

\section{The Mourre theory}\label{s:mourre}
In \cite{Pu}, C.R.\ Putnam used a positive commutator estimate to insure that the spectrum of an operator is purely absolutely continuous. His method was unfortunately not very flexible and did not allow the presence of eigenvalue. In  \cite{Mo81,Mo83}, E.\ Mourre had the idea to localise in energy the positive commutator estimate. Thanks to some hypothesis of regularity, he proved that the embedded eigenvalues can accumulated only at some thresholds, that the singularly continuous spectrum is empty and also established a limiting absorption principle, away from the eigenvalues and from the thresholds. Many papers have shown the power of Mourre's commutator theory for a  wide class of self-adjoint operators, e.g., \cite{BFS,BCHM,CGH,DJ,FH,GGM1,GG,HUS,JMP,S}. We refer to \cite{ABG} for the optimised theory and to \cite{GJ1, GJ,Ge} for recent developments.

Let us now, briefly recall Mourre's commutator theory. The aim is to establish some
spectral properties of a given (unbounded) self-adjoint operator $H$ acting in some complex and separable Hilbert space $\Hc$ with the help of an external unbounded and self-adjoint operator $\Ac$. Let $\|\cdot\|$ denote the norm of bounded operators on $\Hc$ and $\sigma (H)$ the spectrum of $H$. Recall that the latter is real. We endow $\Dc(H)$, the domain of $H$, with its graph norm
We denote by $R(z):=(H-z)^{-1}$ the resolvent of $H$ in $z$.
Take an other Hilbert space $\Kc$ such that there is a dense and injective embedding from $\Kc$ to $\Hc$, by identifying $\Hc$ with its antidual $\Hc^*$, we have: $\Kc \hookrightarrow \Hc \simeq \Hc^* \hookrightarrow \Kc^*$, with  dense and injective embeddings.

We introduce some regularity classes with respect to $\Ac$ and follow \cite[Chapter 6]{ABG}.
Given $k\in \N$,  we say that $H\in \Cc^k(\Ac)$ if for all $f\in \Hc$, the map $\R\ni t\mapsto e^{\rm i t\Ac}(H+\rm i)^{-1} e^{-\rm i t\Ac}f\in \Hc$ has the usual $\Cc^k$ regularity.
We say that $H\in \Cc^{k,u}(\Ac)$ if the map $\R\ni t\mapsto e^{\rm i t\Ac}(H+\rmi) e^{-\rm i t\Ac}\in\Bc(\Hc)$  is of class  $\Cc^k(\R,\Bc(\Hc))$, where $\Bc(\Hc)$ is endowed with the norm operator topology.

We start with an example, e.g., \cite[Proposition 2.1]{GJ1}.
\begin{lemma}
\label{Lemma:2}
For $\phi, \varphi \in \Dc(A)$, the rank one operator $|\phi \rangle \langle \varphi | : \psi \mapsto \langle \varphi, \psi \rangle \phi$ is of class $\Cc^{1}(A)$ and
\[\left[|\phi \rangle \langle \varphi |, A\right]= |\phi \rangle \langle A \varphi | - |A \phi \rangle \langle \varphi |.\]
By induction, given $n\in \N$ and $\phi, \varphi \in \Dc(A^n)$,
$|\phi \rangle \langle \varphi |$ is of class $\Cc^{n}(A).$
\end{lemma}
We turn to a criterion in term of commutator.

 \begin{theorem}[{\cite[p.258]{ABG}}]  Let $\Ac$ and $H$ be two self-adjoint operators in the Hilbert space $\Hc$.
 The following points are equivalent: \begin{enumerate}
\item $H\in\Cc^1(\Ac)$.
\item\label{3.1} For one (then for all) $z\not\in\sigma(H)$, there is a finite $c$ such that
\begin{align*}|\langle\Ac f,R(z)f\rangle-\langle R(\overline{z})f&,\Ac f\rangle|\leq c\|f\|^2, \mbox{ for all } f\in\Dc(\Ac).\end{align*}
\item
\begin{enumerate}
\item\label{3.2}  There is a finite $c$ such that for all $f\in\Dc(\Ac)\cap\Dc(H)$:
\begin{align*}|\langle\Ac f,Hf\rangle-\langle Hf,\Ac f\rangle|&\leq c(\|Hf\|^2+\|f\|^2).\end{align*}
\item  For some (then for all) $z\not\in\sigma(H)$, the set \[\{f\in\Dc(\Ac), R(z)f\in\Dc(\Ac) \mbox{ and } R(\overline{z})f\in\Dc(\Ac)\}\mbox{is a core for } \Ac.\]
\end{enumerate}
\end{enumerate}
\end{theorem}
Note that \eqref{3.1} yields that the commutator $[\Ac, R(z)]$ extends to a bounded operator in the form
sense. We shall denote the extension by $[\Ac, R(z)]_\circ$. In the same
way, from \eqref{3.2}, the commutator $[H, \Ac]$ extends to a unique
element of $\Bc\big(\Dc(H), \Dc(H)^*\big)$ denoted by $[H,\Ac]_\circ$. Note that $\Dc(H)$ is endowed with the graph norm of $H$ and that $\Dc(H)^*$ denotes its anti-dual. Moreover, if $H\in \Cc^1(\Ac)$ and $z\notin \sigma(H)$,
\begin{align*}
\big[A, (H-z)^{-1}\big]_\circ &=\quad  \underbrace{(H-z)^{-1}}_{\Hc
  \leftarrow \Dc(H)^*}\quad  \underbrace{[H, A]_\circ}_{\Dc(H)^*\leftarrow
  \Dc(H)} \quad \underbrace{(H-z)^{-1}}_{\Dc(H)\leftarrow \Hc}.
\end{align*}
Here, we use the Riesz lemma to identify $\Hc$ with its anti-dual
$\Hc^*$.

Note that, in practice, the condition (3.b) could be delicate to check. This is addressed by the next lemma.
\begin{lemma}[{\cite[Lemma A.2]{GoMo}}]\label{l:GoMo}
Let $\Dr$ be a subspace of $\Hc$ such that $\Dr\subset\Dc(H)\cap\Dc(\Ac)$, $\Dr$ is a core for A and
$H\Dr\subset\Dr$. Let $(\chi_n)_{n\in\N}$ be a family of bounded operators such that
\begin{enumerate}
\item $\chi_n\Dr\subset\Dr$, $\chi_n$ tends strongly to $1$ as $n\to\infty$, and $\sup_n\|\chi_n\|_{\Bc(\Dc(H))}<\infty$.
\item $\Ac\chi_n f\to \Ac f$, for all $f\in\Dr$, as $n\to\infty$.
\item There is $z\not\in\sigma(H)$, such that $\chi_n R(z)\Dr\subset\Dr$ and $\chi_n R(\overline{z})\Dr\subset\Dr$.
\end{enumerate}
Suppose also that for all $f\in\Dr$
\[\lim_{n\to\infty}\Ac[H,\chi_n]R(z)f=0  \mbox{ and } \lim_{n\to\infty}\Ac[H,\chi_n]R(\overline{z})f=0 .\]
Finally, suppose that there is a finite $c$ such that
\[|\langle\Ac f, Hf\rangle-\langle Hf,\Ac f\rangle|\leq c(\|Hf\|^2+\|f\|^2), \mbox{ for all } f\in\Dr.\]
Then, one has $H\in\Cc^1(\Ac)$.
\end{lemma}

 We define other refined classes of regularity: 
 \begin{align*}
\mbox{We say that } H\in\Cc^{0,1}(\Ac)& \mbox{ if }\int^1_0 \left\| [(H+\rm i)^{-1},e^{\rm i t\Ac}]\right\| \frac{dt}{t}<\infty.
 \\
\mbox{We say that } H\in\Cc^{1,1}(\Ac)& \mbox{ if }\int^1_0 \left\| [[(H+\rm i)^{-1},e^{\rm i t\Ac}],e^{\rm i t\Ac}]\right\| \frac{dt}{t^2}<\infty.
\end{align*}
Thanks to \cite[p. 205]{ABG}, it turns out that

Given  an interval open interval $\Ic$, we denote by $E_{\Ic}(H)$ the spectral projection of $H$ above $\Ic$. We say that the \emph{Mourre estimate} holds true for $H$ on $\Ic$ if there exist $c>0$ and a compact operator $K$ such that
\begin{align}\label{eq:mourre}
E_\Ic (H)[H, \rmi \Ac]_{\circ}E_\Ic(H)\geq E_\Ic(H)\, (c\, +\, K)\, E_\Ic(H),
\end{align}
when the inequality is understood in the form sense.
We say that we have a \emph{strict Mourre estimate} holds for $H$ on the open interval $\Ic'$ when there exists $c'>0$
such that
\begin{align}\label{eq:mourrestrict}
E_{\Ic'} (H)[H, \rm i \Ac]_{\circ}E_{\Ic'}(H)\geq c' E_{\Ic'}(H).
\end{align}
Assuming $H\in \Cc^1(\Ac)$,  \eqref{eq:mourre}, and $\lambda\in \Ic$ is not an eigenvalue, therefore there exists an open interval $\Ic'$ that contains $\lambda$ and $c'>0$ such that \eqref{eq:mourrestrict}.
The aim of Mourre's commutator theory is to
show a \emph{limiting absorption principle} (LAP), see \cite[Theorem 7.6.8]{ABG}.
\begin{theorem}Let $H$ be a self-adjoint operator, with $\sigma(H)\neq \R$. Assume that $H\in\Cc^1(\Ac)$ and the Mourre estimate  \eqref{eq:mourre} holds true for $H$ on $\Ic$. Then
\begin{enumerate}
\item The number of eigenvalues (counted with multiplicity) of $H$, that are in $\Ic$, is finite.
\end{enumerate}
Assuming furthermore that $K=0$ in \eqref{eq:mourre}, it yields:
\begin{enumerate}
\item[(2)] $H$ has no eigenvalues in $\Ic$.
\item[(3)]  If $H\in\Cc^{1,1}(\Ac)$ and $K=0$, $s>1/2$ and $\Ic'$ a compact sub-interval of $\Ic$, then \[\sup_{\Re(z)\in\Ic',\Im(z)\neq0}\|\langle \Ac\rangle^{-s}(H-z)^{-1}\langle \Ac\rangle^{-s}\| \mbox{ exists and finite}.\]
\end{enumerate}
Moreover, in the norm topology of bounded operators, the boundary values of the resolvent:
\[\Ic'\ni\lambda\mapsto\lim_{\rho\to0^{\pm}}\langle \Ac\rangle^{-s}(H-\lambda-\rmi\rho)^{-1}\langle \Ac\rangle^{-s} \mbox{ exists and continuous}.\]
\end{theorem}
For more details and deeper results, see \cite[Proposition 7.2.10, Corollary 7.2.11, Theorem 7.5.2]{ABG}.

\section{The free model}
\subsection{Construction of the graph}
\label{s:setting}
We discuss two different product of graphs. To start off,
given $\Gc_1:=(\Ec_1, \Vc_1, m_1)$ and
$\Gc_2:=(\Ec_2, \Vc_2, m_2)$, the \emph{Cartesian product of $\Gc_1$ by $\Gc_2$} is defined by
$\Gc^{\diamond}:=(\Ec^{\diamond}, \Vc^{\diamond}, m^{\diamond})$, where $\Vc^{\diamond}:= \Vc_1\times \Vc_2$,
\begin{align*}
\left\{\begin{array}{rl}
m^{\diamond}(x,y):=& m_1(x)\times m_2(y),
\\
\Ec^{\diamond}\left((x,y),(x',y')\right):=
&\Ec_1(x,x')\times \delta_{y, y'}  m_2(y)+
m_1(x)\delta_{x,x' } \times \Ec_2(y,y').
\end{array}\right.
\end{align*}
We denote it by $ \Gc_1\times \Gc_2:= \Gc^{\diamond}$. This definition generalises the unweighted Cartesian product, e.g., \cite{Ha}. It is used in several places in the literature, e.g., see \cite[Section 2.6]{Ch} and see  \cite{BGKLM} for a generalisation.
\\ The terminology is motivated by the following decomposition:
\begin{align*}
\Delta_{\Gc^{\diamond}} = \Delta_{\Gc_1} \otimes 1 +
1 \otimes \Delta_{\Gc_2},
\end{align*}
where $\ell^2(\Vc,m)\simeq\ell^2(\Vc_1, m_1)\otimes \ell^2(\Vc_2, m_2)$. Note that
\[e^{\rmi t \Delta_{\Gc^{\diamond}}} = e^{\rmi t \Delta_{\Gc_1}} \otimes e^{\rmi t \Delta_{\Gc_2}}, \quad \forall t\in \R.\]
 We refer to \cite[Section VIII.10]{RS} for an introduction to
 the tensor product of self-adjoint operators.

We now introduce a \emph{twisted Cartesian product}. We refer to
 \cite[Section 2.2]{GT} for motivations, its link with hyperbolic geometry and generalisations.
Given $\Gc_1:=(\Ec_1, \Vc_1, m_1)$ and
$\Gc_2:=(\Ec_2, \Vc_2, m_2)$, we define the \emph{product of $\Gc_1$ by $\Gc_2$} by
$\Gc:=(\Ec, \Vc, m)$, where $\Vc:= \Vc_1\times \Vc_2$ and
\begin{align*}
\left\{\begin{array}{rl}
m(x,y):=& m_1(x)\times m_2(y),
\\
\Ec\left((x,y),(x',y')\right):=
&\Ec_1(x,x')\times \delta_{y, y'} +
\delta_{x,x' } \times \Ec_2(y,y'),
\end{array}\right.
\end{align*}
for all $x, x'\in \Vc_1$ and $y, y'\in \Vc_2$. We denote $\Gc$ by
$\Gc_1\times_{\tau} \Gc_2$. If $m=1$,  note that $\Gc_1\times_{\tau} \Gc_2 = \Gc_1\times \Gc_2$.

Under the representation $\ell^2(\Vc, m)\simeq \ell^2(\Vc_1, m_1)\otimes
\ell^2(\Vc_2, m_2)$,
\begin{align}\label{e:deg_cp}
\deg_{\Gc_1\times_{\tau} \Gc_2}(\cdot)= \deg_{\Gc_1}(\cdot)\otimes \frac{1}{m_2(\cdot)}
+ \frac{1}{m_1(\cdot)}\otimes \deg_{\Gc_2}(\cdot)
\end{align}
and
\begin{align}\label{e:rule}
\Delta_{\Gc_1\times_{\tau} \Gc_2} = \Delta_{\Gc_1} \otimes \frac{1}{m_2(\cdot)} +
\frac{1}{m_1(\cdot)} \otimes \Delta_{\Gc_2}.
\end{align}
If $m$ is non-trivial, we stress that the Laplacian obtained with our
product is usually not unitarily equivalent to the Laplacian obtained
with the Cartesian product.

A hyperbolic manifold of finite volume is the union of a compact part, of a cusp, and a funnel, e.g., \cite[Theorem 4.5.7]{Th}. In this article we study a discrete analog.
\vspace*{0.5cm}
\begin{center}
\bf{In the sequel, we take $m_2$ constant on $\Vc_2$.}
\end{center}
\vspace*{0.5cm}
The graph $\Gc:= (\Ec, \Vc, m)$ is divided into three parts: A \emph{cusp} part, a \emph{funnel} part, and a finite part.
Set $\Gc^\star:= (\Ec^\star, \Vc^\star, m^\star)$ be the induced graph of $\Gc$ over $\Vc^\star$ where $\star \in \{{\rm c, f, 0}\}$ and $\Vc$ is a disjoint reunion of $\Vc^{\rm c}$, $\Vc^{\rm f}$, and $\Vc^{\rm 0}$.

We consider $\Gc_1^{\rm c}:=(\Ec_1^{\rm c}, \Vc_1^{\rm c}, m_1^{\rm c})$, where
\[\Vc_1^{\rm c}:=\N, \quad m_1^{\rm c}(n):=\exp(-n), \mbox{ and } \Ec_1^{\rm c}(n, n+1):= \exp(-(2n+1)/2),\]
for all $n\in \N$ and $\Gc_2^{\rm c}:=(\Ec_2^{\rm c}, \Vc_2^{\rm c},  m_2)$ a possibly disconnected connected finite graph. 
Set $\Gc^{\rm c}:=\Gc_1^{\rm c}\times_{\tau} \Gc_2^{\rm c}$. This is a \emph{cusp part}. Note it is of finite volume as:
\[\sum_{(x,y)\in \Vc_1^c\times \Vc_2^c} m_{\Gc^{\rm c}} (x,y)<\infty.\]

We consider $\Gc_1^{\rm f}:=(\Ec_1^{\rm f}, \Vc_1^{\rm f}, m_1^{\rm f})$, where
\[\Vc_1^{\rm f}:=\N, \quad m_1^{\rm f}(n):=\exp(n), \mbox{ and } \Ec_1^{\rm f}(n, n+1):= \exp((2n+1)/2),\]
for all $n\in \N$ and $\Gc_2^{\rm f}:=(\Ec_2^{\rm f}, \Vc_2^{\rm f}, m_2)$ a connected finite graph. 
 Set $\Gc^{\rm f}:=\Gc_1^{\rm f}\times_{\tau} \Gc_2^{\rm f}$. This is a \emph{funnel part}.

For the \emph{compact part}, we ask that for all $x\in \Vc^{0}$, $\supp(\Ec(x, \cdot))$ is finite and $m^{\rm 0}(x)>0$.

We now, take advantage of \[\ell^2(\Gc):=\ell^2(\Gc^{\rm f})\oplus\ell^2(\Gc^{\rm 0})\oplus\ell^2(\Gc^{\rm c}).\] We have that \[\Delta_{\Gc}:=\Delta_{\Gc^{\rm f}}\oplus\rm0\oplus\Delta_{\Gc^{\rm c}}+K_0,\]
where $K_0$ is an operator of finite rank with support in $\Cc^c(\Vc)$.

To analyse the perturbations of operator we shall rely on the following gauge transformation, e.g., \cite{Go, CTT, HK}. See also \cite{BG} for some historical references.
\begin{proposition}\label{p:uni}
Let $\Gc:=(\Vc,  \Ec, m)$ be a weighted graph and $m:\Vc\to(0,
\infty)$ be a weight. The following map is unitary:
\begin{align}
\nonumber
T_{m\rightarrow m'}f:\ell^2(\Vc, m)&\to \ell^2(\Vc, m')
\\
\label{e:T}
f&\mapsto \left(x \mapsto \sqrt{\frac{m(x)}{m'(x)}} f(x)\right).
\end{align}
We have:
\begin{align}\label{e:transform}
\Delta_{\Gc'}^{\Fc}=T_{m\rightarrow m'}\left(\Delta_{\widetilde{\Gc}}- W(\cdot)\right)^{\Fc} T_{m\rightarrow m'}^{-1},
\end{align}
where $\Gc':=(\Vc, \Ec', m')$, $\widetilde{\Gc}:=(\Vc, \widetilde{\Ec}, m)$ and,
\begin{align*}
\widetilde{\Ec}(x,y)&:= \Ec' (x,y) \sqrt{\frac{m(x) m(y)}{m'(x) m'(y)}}
\\
W(x)&:=\frac{1}{m(x)}\sum_{y\in \Vc} \widetilde{\Ec}(x,y)  \left(1
 - \sqrt{\frac{m(x)m'(y)}{m(y)m'(x)}}\right).
\end{align*}
Here we emphasised the choice Friedrichs extension with the symbol $\Fc$.
\end{proposition}
\subsection{Mourre estimate on $\N$ }\label{Mourre N}
In this section we make a preliminary work on the half axis. We construct a conjugate operator, prove a Mourre estimate for $\Delta_\N$ and check the regularity conditions. This is a known result, e.g., \cite{AF}, see also \cite{GG, Mic}.

Given $f\in \ell^{2}(\N,1),$ we set
\[\forall n\in\N^*, \quad Uf(n):=f(n-1) \text{ and } Uf(0):=0.\]
Note that $U^* f(n)=f(n+1), \forall n\in\N$. The operator $U$ is an isometry and is not unitary: we have $U^{*}U=\id$ and $UU^{*}=\un_{[1, \infty[}(\cdot)$.

We define by $Q$ the operator of multiplication by $n$ in $\ell^2(\N,1)$.
Namely, it is the closure of the operator given by  $(Qf)(n)= nf(n)$ for all $n\in\N$ and $f\in \Cc_{c}(\N)$.
It is essentially self-adjoint on $\Cc_{c}(\N)$.
In \cite{GG}, one finds the following elementary relations:
\begin{align}\label{dec}
QU = U(Q+1), \ U^*Q = (Q+1) U^*  \mbox{ and } UQU=U^2(Q+1) \quad \mbox{on } \Dc(Q).\end{align}
The operator $\Delta_{\N}$ is defined by \eqref{e:DeltaG},
 where $\N\simeq(\N,\Ec_{\N},m)$, with $\Ec_{\N}(n,n+1)=1$ and $m(n)=1$ for all $n\in\N$.
Explicity, we have
\[\Delta_\N f(n):=
\left\{\begin{array}{cl}
2f(n)-f(n-1)-f(n+1)& \text{ if } n\geq 1,
\\
f(n)- f(n+1)& \text{ if } n=0,
\end{array}\right.
\quad \forall f\in \ell^2(\N, 1).\]
We can express it with the help of $U$. Namely, we have:
\[\Delta_{\N} = 2 - (U+U^*) - \un_{\{0\}}(\cdot).\]
A standard result is :
\[\sigma_{\rm ess}(\Delta_\N)= [0,4] \quad \mbox{ and } \quad \sigma_{\rm sc}(\Delta_\N)= \emptyset.\]
We construct the conjugate operator in $\ell^2(\N, 1)$. On the space $\Cc_c(\N)$, we define
\begin{align}\nonumber
  \Ac_\N|_{\Cc_c(\N)}&:=\frac{1}{2} \left( SQ +  QS\right), \quad \mbox{ where } S:= \frac{U-U^*}{2\rmi}
  \\ \nonumber
  &=\frac{\rmi}{2}\left(U\left(Q+\frac{1}{2}\right)-U^*\left(Q-\frac{1}{2}\right)\right)
 \\ \label{e:AN}
&=-\frac{\rmi}{2}\left(\frac{1}{2}\left(U^{*}+U\right)+Q\left(U^{*}-U\right)\right).
\end{align}
We denote by $\Ac_\N$ its closure.
\begin{lemma}\label{l:ANessaa}
The operator $\Ac_\N$ is essentially self-adjoint on $\Cc_c(\N)$ and
\[\Dc(\Ac_\N)=\Dc(QS):= \{f\in \ell^2(\N), Sf\in \Dc(Q)\}.\]
\end{lemma}
We refer to \cite{GG} and \cite[Lemma 5.7]{Mic} for the essential self-adjointness and \cite[Lemma 3.1]{GG} for the domain.

We give a first technical lemma.
\begin{lemma}\label{l:A0}
On $\Cc_c(\mathbb{N})$, we have
\begin{align*}
(U^{*}+U)\Ac_{\mathbb{N}}&=-\frac{\rmi}{2}\left((U^{2*}-U^{2})Q-\un-\frac{1}{2}(\un_{\{0\}}(\cdot)+U^{2}+U^{*2})\right), \\
\Ac_{\mathbb{N}}(U^{*}+U)&=\frac{\rmi}{2}\left((U^{2}-U^{2*})Q-\frac{1}{2}(U^{2}+U^{*2})-\un-\frac{1}{2}\un_{\{0\}}(\cdot)\right).
\end{align*}
\end{lemma}
\proof We compute on $\Cc_c(\N)$. The statement follows easily from
\begin{align*}  U\Ac_{\mathbb{N}}
&=-\frac{\rmi}{2}\left((\un_{\mathbb{N}}-U^{2})Q-\frac{1}{2}(\un_{\mathbb{N}}+U^{2}),\right).
\\
U^{*}\Ac_{\mathbb{N}}&=-\frac{\rmi}{2}\left((U^{*}-\un)Q-\frac{1}{2}(U^{*2}-\frac{1}{2})\right).\end{align*}
by taking the adjoint. \qed

We can compute the first comutator.
\begin{lemma}\label{l:commuA}
The operator $\Delta_\N$ is $\Cc^1(\Ac_\N)$ and we have:
\begin{equation}\label{e:commuA}
[\Delta_{\mathbb{N}},\rmi \Ac_{\mathbb{N}}]_\circ=\frac{1}{2}\Delta_{\mathbb{N}}(4-\Delta_{\mathbb{N}}) + K_1,
\end{equation}
with $K_1$ a finite rank operator belonging to $\Cc^\infty(A)$.
\end{lemma}
This lemma is essentially given in \cite{GG}, see also \cite{AF} for another type of presentation. For the convenience of the reader we reproduce it.
\proof First, since $\delta_{\{0\}}\in \Dc(\Ac^n)$ for all $n\in \N$, $\delta_{\{0\}}$ and $K_1:=[\delta_{\{0\}},\rm i\Ac_{\mathbb{N}}]_\circ$ belong to $\Cc^1(\Ac_\N)$ by Lemma \ref{Lemma:2}.
Next, we turn to the other part and work in the form sense and by density. Let $f\in \Cc_c(\N)$. Since $\Delta_\N f \in \Cc_c(\N)$ and using Lemma \ref{l:A0}, we obtain:
 \begin{align*}
 \langle f,[\Delta_{\mathbb{N}},\rmi \Ac_{\mathbb{N}}] f\rangle&:=
 \langle \Delta_\N f, \rmi \Ac_\N f\rangle - \langle -\rmi \Ac_\N f, \Delta_\N f \rangle
\\
&=\rmi\langle f,\Ac_\N(U^{*}+U)-(U^{*}+U)\Ac_{\mathbb{N}} f \rangle +\langle f,[\delta_{\{0\}},\rmi \Ac_{\mathbb{N}}] f \rangle
\\
&= \frac{1}{2}\langle f, \Delta_{\mathbb{N}}(4-\Delta_{\mathbb{N}}) f\rangle +\langle f,[\delta_{\{0\}},\rmi \Ac_{\mathbb{N}}]_\circ f \rangle.
\end{align*}
Since $\Delta_{\mathbb{N}}(4-\Delta_{\mathbb{N}})$ and $[\delta_{\{0\}},\rm i\Ac_{\mathbb{N}}]_\circ$
 are bounded operators and since $ \Cc_c(\N)$ is a core for $\Ac_\N$, there is a constant $c$ such that
 \[|\langle \Delta_\N f, \rmi \Ac_\N f\rangle - \langle -\rmi \Ac_\N f, \Delta_\N f \rangle |\leq c\| f\|^2, \mbox{ for all } f\in \Dc(\Ac).\]
 Hence, it is $\Cc^1(\Ac_\N)$. By density, we also obtain \eqref{e:commuA}.\qed

By induction, we infer:
\begin{corollary} $\Delta_{\N} \in \Cc^\infty(\Ac_\N)$.
\end{corollary}
We mention \cite{Mic} for an anisotropic use on $\Z$ based on the Mourre theory of $\Delta_\N$.

\subsection{The funnel side}\label{section3}In this section we construct a conjugate operator for $\Delta_{\Gc^{\rm f}}$ and establish a Mourre estimate.
\subsubsection{A first step into the analysis}
As seen above, under the identification
\begin{align}\label{e:otimes}
\ell^{2}(\Vc^{\rm f},m)&= \ell^{2}(\N,m^{\rm f}_1)\otimes  \ell^{2}(\Vc^{\rm f}_2, m_2^{\rm f}).
\end{align}
We have
\begin{align}\label{e:deltaf}
\Delta_{\Gc^{\rm f}}:= \Delta_{\Gc_1^{\rm f}} \otimes \frac{1}{m_2^{\rm f}}+
\frac{1}{m_1^{\rm f}(\cdot)} \otimes  \Delta_{\Gc_2^{\rm f}}.
\end{align}
Recall here that $m_2$ is a constant.
The first remark is that
\begin{lemma}\label{l:compf}
\[\frac{1}{m_1^{\rm f}(\cdot)} \otimes  \Delta_{\Gc_2^{\rm f}}\in \Kc\left(\ell^2(\Vc^{\rm f})\right).\]
\end{lemma}
\proof Note that $\Delta_{\Gc_2^{\rm f}}$ is of finite rank since $\Vc_2$ is finite and that $\displaystyle\frac{1}{m_1^{\rm f}(\cdot)}$ is a compact operator since $m_1^{\rm f}(n)\to \infty$, as $n\to\infty$.\qed

Since $m_2$ is constant and $\deg_{\Gc_1^{\rm f}}$ is bounded, we obtain:
 \begin{proposition}
 We have $\Delta_{\Gc^{\rm f}}\in \Bc\left(\ell^2(\Vc^{\rm f}), m^{\rm f}\right)$.
 \end{proposition}



Recalling the Proposition \ref{p:uni}, we obtain:
\begin{align*}
T_{1\rightarrow m_1^{\rm f}}^{-1}\Delta_{\Gc^{\rm f}_1}T_{1\rightarrow m_1^{\rm f}}=\Delta_{\N}+ (e^{1/2}-1)\un_{\{0\}}+
e^{1/2}+e^{-1/2}-2
\end{align*}
Recalling Lemma \ref{l:compf}, we infer immediately
\[\sigma_{\rm ess}(\Delta_{\Gc^{\rm f}})= \left[\frac{\alpha}{{m_2}}, \frac{\beta}{{m_2}}\right] \quad \text{ and } \quad \sigma_{\rm sc}(\Delta_{\Gc^{\rm f}_1})= \emptyset,\]
with $\alpha$ and $\beta$ are given in \eqref{e:c1c2}.
\subsubsection{Construction of the conjugate operator}
In order to get also $\sigma_{\rm sc}(\Delta_{\Gc^{\rm f}})= \emptyset$, we rely on the Mourre theory and construct a conjugate operator for $\Delta_{\Gc^{\rm f}}$. Recalling \eqref{e:AN} and with respect to \eqref{e:otimes}, we set

\begin{align}\label{Afunnel}
\Ac_{\Gc^{\rm f}}:= 
\Ac_{m_1^{\rm f}}\otimes 1_{\Vc_2^{\rm f}}:=
T_{1\rightarrow m_1^{\rm f}}\Ac_{\N}T^{-1}_{1\rightarrow m_1^{\rm f}}\otimes 1_{\Vc_2^{\rm f}}.
\end{align}

It is essentially self-adjoint on $\Cc_c(\Vc^{\rm f})$ and on $\Cc_c(\N)\otimes \ell^2(\Vc^{\rm f})$ by Lemma \ref{l:ANessaa}. It acts as follows:

\begin{proposition}\label{p:Am1f} On $\Cc_c(\N)$, we have
 \[\Ac_{m_1^{\rm f}}=\frac{\rmi}{2}\left(e^{1/2}(Q-1/2)U-e^{-1/2}(Q+1/2)U^*\right).\]
\end{proposition}
\proof Let $f\in\Cc_c(\N),$
\begin{align*}&\Ac_{m_1^{\rm f}}f(n)
=-\frac{\rmi}{2\sqrt{m_1^{\rm f}(n)}} \left(\frac{1}{2} \left(U+U^*\right) + Q\left(U^*-U\right)\right)T^{-1}_{1\rightarrow m_1^{\rm f}}f(n)
\\
&=\frac{\rmi}{2}\left(\left(n-\frac{1}{2}\right)\sqrt{\frac{m_1^{\rm f}(n-1)}{m_1^{\rm f}(n)}}f(n-1)-\left(n+\frac{1}{2}\right)\sqrt{\frac{m_1^{\rm f}(n+1)}{m_1^{\rm f}(n)}}f(n+1)\right)
\\
&=\frac{\rmi}{2}\left(e^{1/2}\left(n-\frac{1}{2}\right)Uf(n)-e^{-1/2}\left(n+\frac{1}{2}\right)U^*f(n)\right).
\end{align*}
This concludes the proof. \qed

We turn to the regularity. In order to lighten the computation, given a graph $\Gc=(\Ec, \Vc, m)$, we write
\[T_1 \simeq T_2 \text{ if there is  $K:\Cc_c(\Vc) \to \Cc_c(\Vc)$ of finite rank such that } T_1 = T_2 + K \]
Thanks to Lemma \ref{Lemma:2} and Proposition \ref{p:Am1f}, we obtain immediately:
\begin{lemma}\label{l:class}
Assume that  $T_1\simeq T_2$. Then for all $n\in \N$,
\[T_1\in\Dc(\Ac_{\Gc^{\rm f}}^n) \Leftrightarrow T_2\in\Dc(\Ac_{\Gc^{\rm f}}^n).\]
\end{lemma}
 We have:
\begin{lemma}\label{l:C1Af}
We have $\Delta_{\Gc^{\rm f}}\in \Cc^1(\Ac_{\Gc^{\rm f}})$ and
\begin{align}\label{e:C1Af}
 \left[\Delta_{\Gc^{\rm f}}, \rmi \Ac_{\Gc^{\rm f}}\right]_\circ &
= w(\Delta_{\Gc^{\rm f}}) + K,
\end{align}
where
\[w^{\rm f}(x):= \frac{m_2}{2}\left(x-\frac{\alpha}{m_2}\right)\left(\frac{\beta}{m_2}-x\right),\]
with $\alpha$ and $\beta$ as in \eqref{e:c1c2} and $K$ is a compact operator.
\end{lemma}
\proof We prove that $[\Delta_{\Gc^{\rm f}}, \rmi \Ac_{\Gc^{\rm f}}]_\circ\in \Bc(\ell^2(\Vc^{\rm f}, m^{\rm f}))$.
As in Lemma \ref{l:commuA} and working in the form sense on $\Cc_c(\N)\otimes \ell^2(\Vc_2^{\rm f})$, a straightforward computation leads to
\begin{align}
\nonumber
\left[\Delta_{\Gc^{\rm f}_1}\otimes \frac{1}{m_2}, \rmi \Ac_{\Gc^{\rm f}}\right]&\simeq\frac{1}{2} (\Delta_{\Gc^{\rm f}_1}-\alpha) (\beta - \Delta_{\Gc^{\rm f}_1}) \otimes\frac{1}{m_2}
\\
\nonumber
&\hspace*{-2cm}\simeq w^{\rm f}(\Delta_{\Gc^{\rm f}}) - \frac{m_2}{2}\left(\frac{1}{m_1(\cdot)}\otimes\Delta_{\Gc_2^{\rm f}}\right)\left(\frac{\beta}{m_2}-\Delta_{\Gc_1^{\rm f}}\otimes\frac{1}{m_2}-\frac{1}{m_1(\cdot)}\otimes\Delta_{\Gc_2^{\rm f}}\right)
\\
\label{e:C1Afa}
&\hspace*{-2cm}\quad+\frac{m_2}{2}\left(\Delta_{\Gc_1^{\rm f}}\otimes\frac{1}{m_2}-\frac{\alpha}{m_2}\right)\left(\frac{1}{m_1(\cdot)}\otimes\Delta_{\Gc_2^{\rm f}}\right).
\\
\nonumber
&\hspace*{-2cm}= w^{\rm f}(\Delta_{\Gc^{\rm f}})  +K'
\end{align}
where $K'$ is a compact operator coming from Lemma \ref{l:compf} and Lemma \ref{Lemma:2}.
We turn to the second part of $\Delta_{\Gc^{\rm f}}$.
\begin{align}
\nonumber
\left[\frac{1}{m_1^{\rm f}(\cdot)} , \rmi \Ac_{m_1^{\rm f}}\right] \otimes \Delta_{\Gc_2^{\rm f}}&=
T_{1\to m_1^{\rm f}}\left[\frac{1}{m_1^{\rm f}(\cdot)} , \rmi \Ac_\N\right] T_{1\to m_1^{\rm f}}^{-1} \otimes \Delta_{\Gc_2^{\rm f}}
\\
\label{e:C1Afb}
&= T_{1\to m_1^{\rm f}}\left(\frac{1}{2}(e-1) e^{-Q}\left(Q-\frac{1}{2}\right)U\right) T_{1\to m_1^{\rm f}}^{-1} \otimes \Delta_{\Gc_2^{\rm f}},
\end{align}
in the form sense on $\Cc_c(\N)\otimes \ell^2(\Vc_2^{\rm f})$. The operator is a compact  since $U$ is bounded and $\lim_{n\to \infty} e^{-n}(n-1/2)=0$.

This implies that $[\Delta_{\Gc^{\rm f}}, \rmi \Ac_{\Gc^{\rm f}}]_\circ\in \Bc(\ell^2(\Vc^{\rm f}, m^{\rm f}))$ and that \eqref{e:C1Af} holds true.
Finally,  since $\Cc_c(\N)\otimes\ell^2(\Vc_2^{\rm f})$ is a core for $\Ac_{\Gc^{\rm f}}$, we deduce that $\Delta_{\Gc^{\rm f}}\in \Cc^1(\Ac_{\Gc^{\rm f}})$. \qed

\begin{lemma}\label{l:C2Af}
We have $\Delta_{\Gc^{\rm f}}\in \Cc^2(\Ac_{\Gc^{\rm f}})$.
\end{lemma}
\proof As above, since $\Cc_c(\N)\otimes\ell^2(\Vc_2^{\rm f})$ is a core for $\Ac_{\Gc^{\rm f}}$ it is enough to prove that $[[\Delta_{\Gc^{\rm f}}, \rmi \Ac_{\Gc^{\rm f}}]_\circ, \rmi \Ac_{\Gc^{\rm f}}]$, defined initially in the form sense on $\Cc_c(\N)\otimes \ell^2(\Vc_2^{\rm f})$,  extends to an element of $\Bc(\ell^2(\Vc^{\rm f}, m^{\rm f}))$.

We prove that the right hand side of \eqref{e:C1Af} belongs to $\Cc^1(A_{\Gc_{\rm f}})$.
It composed of $w(\Delta_{\Gc_{\rm f}})$ which is $\Cc^1(A_{\Gc_{\rm f}})$ (as product of bounded operators belonging to $\Cc^1(A_{\Gc_{\rm f}})$), terms with finite support that are also in $\Cc^1(A_{\Gc_{\rm f}})$ by Lemma \ref{Lemma:2} and terms similar to \eqref{e:C1Afb}. Therefore $\left[\left[\Delta_{\Gc^{\rm f}_1}\otimes \frac{1}{m_2}, \rmi \Ac_{\Gc^{\rm f}}\right]_\circ, \rmi \Ac_{\Gc^{\rm f}}\right]$ extends to a bounded operator.

We turn to the second part. It remains to show that the left hand side of \eqref{e:C1Afb} belongs to $\Cc^1(A_{\Gc_{\rm f}})$. Repeating the computation done in \eqref{e:C1Afb}, we see that since $\lim_{n\to \infty} e^{-n}\langle n\rangle^2=0$,
$\left[\left[\frac{1}{m_1^{\rm f}(\cdot)} , \rmi \Ac_{m_1^{\rm f}}\right]_\circ,  \Ac_{m_1^{\rm f}}\right]$ extends to a compact operator. \qed

\begin{remark} By induction, we can prove that  $\Delta_{\Gc^{\rm f}}\in \Cc^\infty(\Ac_{\Gc^{\rm f}})$.
\end{remark}

Finally, we establish the Mourre estimate.
\begin{proposition}We have $\Delta_{\Gc^{\rm f}}\in\Cc^2(\Ac_{\Gc^{\rm f}})$. Given a compact interval $\Ic\subset (\alpha/ m_2,\beta/ m_2)$, there are $c>0$, a compact operator $K$ such that
 \begin{align}
 \label{e:mourref}
 E_\Ic(\Delta_{\Gc^{\rm f}})[\Delta_{\Gc^{\rm f}}, \rmi \Ac_{\Gc^{\rm f}}]_\circ E_\Ic(\Delta_{\Gc^{\rm f}})&\geq c E_\Ic(\Delta_{\Gc^{\rm f}}) + K,\end{align} in the form sense. In particular, $\sigma_{\rm sc}(\Delta_{\Gc^{\rm f}})=\emptyset$.
\end{proposition}
\proof
Lemma \ref{l:C2Af} gives that $\Delta_{\Gc^{\rm f}}\in \Cc^2(\Ac_{\Gc^{\rm f}})$. By \eqref{e:C1Af}, we obtain
\begin{align*}
 E_\Ic(\Delta_{\Gc^{\rm f}})[\Delta_{\Gc^{\rm f}}, \rmi \Ac_{\Gc^{\rm f}}]_\circ E_\Ic(\Delta_{\Gc^{\rm f}})&=
       E_\Ic(\Delta_{\Gc^{\rm f}})w(\Delta_{\Gc^{\rm f}})E_\Ic(\Delta_{\Gc^{\rm f}})+K       \\
       &\geq c E_\Ic(\Delta_{\Gc^{\rm f}}) + K, \end{align*}
where $K$ is a compact operator and
\[c:=\frac{m_2}{2}\inf_{x\in I}\left(x-\frac{\alpha}{m_2}\right)\left(\frac{\beta}{m_2}-x\right)>0.\]
The absence of singular continuous spectrum follows from the general theory.  \qed

To lighten the text we did not expand more consequences of the Mourre theory in this case and refer to Theorem \ref{t:mainth} for them.

\subsection{The cusps side}\label{section4}
In this section we construct a conjugate operator for $\Delta_{\Gc^{\rm c}}$ and establish a Mourre estimate. By contrast with the funnel side, we shall refine the tensor product decomposition.

\subsubsection{The model and the low/high energy decomposition} Again we rely on the decomposition
\begin{align}\label{e:otimesc}
\ell^{2}(\Vc^{\rm c},m)&= \ell^{2}(\N,m^{\rm c}_1)\otimes  \ell^{2}(\Vc^{\rm c}_2, m_2).
\end{align}
We have
\begin{align}\label{e:deltac}
\Delta_{\Gc^{\rm c}}:= \Delta_{\Gc_1^{\rm c}} \otimes \frac{1}{m_2}+
\frac{1}{m_1^{\rm c}(\cdot)} \otimes  \Delta_{\Gc_2^{\rm c}}.
\end{align}
Recall that $m_2$ is a constant.
Unlike with the treatment of $\Delta_{\Gc^{\rm f}}$, we refine the tensor product decomposition. In the spirit of \cite{GoMo, GT}, we denote by  $P^{{\rm
    le}}$  the projection on $\ker(\Delta_{\Gc_2})$ and by $P^{{\rm
    he}}$ is the projection on $\ker(\Delta_{\Gc_2})^\perp$.  Here, \emph{le} stands for \emph{low energy} and
\emph{he} for  \emph{high energy}.
We shall take advantage of
\begin{align}
\nonumber
\ell^{2}(\Vc^{\rm c},m)&:= \Hc^{\rm le} \oplus \Hc^{\rm he}
\\
\label{e:lhdecomp}
&:= \ell^{2}(\N,m_1^{\rm c})\otimes \ker(\Delta_2) \oplus \ell^{2}(\mathbb{N},m_1^{\rm c})\otimes \ker(\Delta_2)^\perp.
\end{align}
The main idea is the continuous spectrum comes from the low energy part of the space whereas the discrete spectrum arises from the high energy part.

We have that $\Delta_{\Gc^{\rm c}}:= \Delta_{\Gc^{\rm c}}^{{\rm
    le}}\oplus \Delta_{\Gc^{\rm c}}^{{\rm he}}$, where
    \begin{align}\label{e:deltale}
\Delta_{\Gc^{\rm c}}^{{\rm le}} := \Delta_{\Gc_1^{\rm c}} \otimes \frac{1}{m_2}P^{{\rm le}},
\end{align}
on $(1 \otimes P^{{\rm le}}) \ell^2(\Vc^{\rm c}, m^{\rm c})$, and
\begin{align}\label{e:deltahe}
\Delta_{\Gc^{\rm c}}^{{\rm he}} := \Delta_{\Gc_1^{\rm c}} \otimes \frac{1}{m_2}P^{{\rm he}}+
\frac{1}{m_1^{\rm c}(\cdot)} \otimes P^{{\rm he}} \Delta_{\Gc_2^{\rm c}},
\end{align}
on $(1 \otimes P^{{\rm he}}) \ell^2(\Vc^{\rm c}, m^{\rm c})$. We stress that $m_2$ is constant.

Unlike $\Delta_{\Gc^{\rm f}}$, $\Delta_{\Gc^{\rm c}}$ is unbounded. More precisely we have:
\begin{proposition}\label{p:cuspessaa}
The operator $\Delta_{\Gc^c}$ is essentially self-adjoint on $\Cc_c(\N)\otimes \ell^2(\Vc_2)$ and on $\Cc_c(\Vc^{\rm c})$. Its domain is given by $\Dc\left(\frac{1}{m_1^{\rm c}(\cdot)}\otimes \Delta_{\Gc_2^c}\right)$.
\end{proposition}
\proof
First $m_1^{\rm c}(\cdot)$ is essentially self-adjoint of $\Cc_c(\N)$. Since $\Delta_{\Gc_2^c}$ is bounded, we infer that  $\frac{1}{m_1^{\rm c}(\cdot)}\otimes \Delta_{\Gc_2^c}$ is essentially self-adjoint on $\Cc_c(\N)\otimes \ell^2(\Vc_2)$. Next, since $\Delta_{\Gc_1^{\rm c}} \otimes \frac{1}{m_2}$ is bounded, $\Delta_{\Gc^c}$ is essentially self-adjoint on $\Cc_c(\N)\otimes \ell^2(\Vc_2)$ by the Kato-Rellich Theorem, e.g., \cite[Theorem X.12]{RS}. The statement with  $\Cc_c(\Vc^{\rm c})$ follows by standard approximations. \qed

Using the notation given in \eqref{e:T}, we see that:
\begin{align*}
T_{m_1^{\rm c}\to 1}\Delta_{\Gc_1^{\rm c}}T^{-1}_{m_1^{\rm c}\to 1}=\Delta_{\N}- (e^{-1/2}-1)\un_{\{0\}}(\cdot) +
e^{1/2}+e^{-1/2}-2 \mbox{ in } \ell^2(\N).
\end{align*}

By using for instance some Jacobi matrices techniques, it is well-known that the  essential spectrum of $\Delta_{\Gc^{\rm c}}^{{\rm le}}$ is purely
 absolutely continuous and 
 \[\sigma_{\rm ac}(\Delta_{\Gc^{\rm c}}^{{\rm le}})= [\alpha, \beta],\]
with multiplicity one, e.g., \cite{Wei}. Recall that $\alpha$ and $\beta$ are defined in \eqref{e:c1c2}.

 We turn to the high energy part.
Using  \cite[Equation (10)]{GT},
\begin{align*}
\frac{1}{m_1^{\rm c}(\cdot)} \otimes \Delta_{\Gc_2^{\rm c}}P^{{\rm
    he}} \leq \Delta_{\Gc^{\rm c}} (1\otimes P^{{\rm he}})\leq 2M+ \frac{1}{m_1^{\rm c}(\cdot)}
\otimes \Delta_{\Gc_2^{\rm c}}P^{{\rm he}}.
\end{align*}
Using the min-max Theorem and since $m_1^{\rm c}(n)\to 0$ as $n\to \infty$,
$\Delta_{\Gc^{\rm c}} (1\otimes P^{{\rm he}})$ has a compact resolvent.
We infer that
\[\sigma_{\rm ac}(\Delta_{\Gc^{\rm c}})= \left[\frac{\alpha}{m_2}, \frac{\beta}{m_2}\right] \text{ and } \sigma_{\rm sc}(\Delta_{\Gc^{\rm c}})=\emptyset.\]

\subsubsection{The conjugate operator}
We pursue the analysis of  $\Delta_{\Gc^{\rm c}}$ in order to apply the Mourre theory to it.
We go back to $\ell^2(\N, m_1^{\rm c})\otimes \ker(\Delta_2)$. We set:
\begin{align} \Ac_{\Gc^{\rm c}}^{\rm le}:= T_{m_1^{\rm c}\to 1}^{-1} \Ac_\N T_{m_1^{\rm c}\to 1} \otimes P^{\rm le}.
\end{align}
It is self-adjoint. Straightforwardly we get
\[ \Ac_{\Gc^{\rm c}}^{\rm le}= -\frac{\rmi}{2} \left(e^{-1/2} \left(Q+\frac{1}{2}\right)U^* +e^{1/2} \left(\frac{1}{2}-Q\right)U  \right)\otimes P^{\rm le}
\]
on $\Cc_c (\N)\otimes \ker(\Delta_2)$.
With respect to \eqref{e:lhdecomp}, we set
\[\Ac_{\Gc^{\rm c}}:= \Ac^{\rm le}_{\Gc^{\rm c}}\oplus \Ac^{\rm he}_{\Gc^{\rm c}},  \text{ where } \Ac^{\rm he}_{\Gc^{\rm c}} := 0.\]
By Lemma \ref{l:ANessaa}, it is essentially self-adjoint on $\Cc_c(\N)\otimes \ell^2(\Vc_2^{\rm c}, m_2)$  and also on $\Cc_c(\Vc^{\rm c})$ by standard approximation.
Keeping the notation of Lemma \ref{l:class}, we obtain:
\begin{lemma}\label{l:C1cusp}
We have $\Delta_{\Gc^{\rm c}}\in \Cc^1(\Ac_{\Gc^{\rm c}})$ and
\begin{align} \label{e:C1cusp}
[\Delta_{\Gc^{\rm c}},\rmi \Ac_{\Gc^{\rm c}}]_\circ
\simeq w^{\rm c}(\Delta_{\Gc^{\rm c}}^{\rm le})\oplus 0,
\end{align}
with respect to \eqref{e:lhdecomp}, with $K\in \Kc(\ell^2(\Vc^{\rm c}, m^{\rm c}))$ and
\[w^{\rm c}(x):= \frac{m_2}{2}\left(x-\frac{\alpha}{m_2}\right)\left(\frac{\beta}{m_2}-x\right).\]
In particular, $[\Delta_{\Gc^{\rm c}},\rmi \Ac_{\Gc^{\rm c}}]_\circ\in \Bc(\ell^{2}(\Vc^{\rm c},m^{\rm c}))$.
\end{lemma}
\proof As in Lemma \ref{l:commuA}, using Lemma \ref{Lemma:2}, and working in the form sense on $\Cc_c(\N)\otimes \ell^2(\Vc_2^{\rm c})$, a straightforward computation leads to
\begin{align}
\nonumber
\left[\Delta_{\Gc^{\rm c}_1}\otimes \frac{1}{m_2^c}, \rmi \Ac_{\Gc^{\rm c}}\right]&\simeq\frac{1}{2} (\Delta_{\Gc^{\rm c}_1}-\alpha) (\beta - \Delta_{\Gc^{\rm c}_1}) \otimes\frac{1}{m_2}P^{\rm le}
\\
\label{e:C1Aca}
&\simeq w^c(\Delta_{\Gc^{\rm f}}^{\rm le})\oplus 0.
\end{align}
We turn to the second part of $\Delta_{\Gc^{\rm c}}$.
\begin{align}
\label{e:C1Acb}
\left[\frac{1}{m_1^{\rm c}(\cdot)}\otimes \Delta_{\Gc_2^{\rm f}} , \rmi \Ac_{\Gc^{\rm c}}\right]&= \left[ \frac{1}{m_1^{\rm c}(\cdot)}, \Ac_{\Gc^{\rm c}_1} \right] \otimes 0 =0
\end{align}
This implies that $[\Delta_{\Gc^{\rm c}}, \rmi \Ac_{\Gc^{\rm c}}]_\circ\in \Bc(\ell^2(\Vc^{\rm c}, m^{\rm c}))$ and \eqref{e:C1cusp}.

It remains to prove that $\Delta_{\Gc^{\rm c}}\in \Cc^1(\Ac_{\Gc^{\rm c}})$. We check the hypotheses of Lemma \ref{l:GoMo}. Let $\{\mathcal{X}_n\}_{n\in\N}$ be a family of functions defined on $\Vc_1^{\rm c}\times \Vc_2^{\rm c}$ as follows:
\[\mathcal{X}_n(x_1,x_2):=\left(\left(1-\frac{x_1-n}{n^{2}+1}\right)\vee0\right)\wedge1.\]
Note that $\supp(\mathcal{X}_n)=\lint0,{n^2+n}\rint\times\Vc_2$ and $\forall (x_1, x_2)\in \lint 0,n\rint \times \Vc_2^{\rm c},~\mathcal{X}_n(x_1,x_2)=1$. We set $\Dr:=\Cc_c(\Vc^{\rm c})$.

1) We have $\|\mathcal{X}_n\|_{\infty}=1$ then $\|\mathcal{X}_n(\cdot)\|_{\Bc(\ell^2(\Vc^{\rm c},m^{\rm c}))}=1$.
Moreover, $\mathcal{X}_n(\cdot)$ tends strongly to $1$ as $n\to +\infty$.
  Now, we shall show that $\sup_n\|\mathcal{X}_n(\cdot)\|_{\Dc(\Delta_{\Gc^{\rm f}})}<\infty$.
Since
\begin{align*}
\left[\Delta_{\Gc_1^{\rm c}}\otimes\frac{1}{m_2}+\frac{1}{m_1^{\rm c}(\cdot)}\otimes\Delta_{\Gc_2^{\rm c}},\mathcal{X}_n(\cdot)\right]&
=\underbrace{[\Delta_{\Gc_1^{\rm c}},\mathcal{X}_n(\cdot)]\otimes\frac{1}{m_2}}_{\hbox{bounded by $2\|\Delta_{\Gc_1^c}\|/ m_2$}}+\underbrace{\left[\frac{1}{m_1^{\rm c}(\cdot)},\mathcal{X}_n(\cdot)\right]}_{\hbox{=0}}\otimes\Delta_2,
\end{align*}
then there is $c>0$ such that, for all $f\in\Cc_c(\Vc^{\rm c})$ such that $f\in (\Delta_{\Gc^{\rm c}}+\rmi)\Cc_c(\Vc^{\rm c})$ and $n\in\N$,
\[\|(\Delta_{\Gc^{\rm c}}+\rmi)\mathcal{X}_n(Q)(\Delta_{\Gc^{\rm c}}+\rmi)^{-1}f\|\leq c \|f\|.\]
Since $\Delta_{\Gc^{\rm c}}$ is essentially self-adjoint on $\Cc_c(\Vc^c)$ and since $-\rmi \notin \sigma(\Delta_{\Gc^{\rm c}})$, it holds for all $f\in\ell^2(\Vc^{\rm c}, m^{\rm c})$. In particular, we derive that
$\|(\Delta_{\Gc^{\rm c}}+\rmi)\mathcal{X}_n(Q)f\|\leq c\|(\Delta_{\Gc^{\rm c}}+\rmi)f\|$,
for all $f\in\ell^2(\Vc^{\rm c},m^{\rm c})$. In particular, $\sup_n\|\mathcal{X}_n(\cdot)\|_{\Dc(\Delta_{\Gc^{\rm f}})}<\infty$.

2)  Given $f\in \Cc_c(\Vc^{\rm c})$, note that for $n$ large enough $\mathcal{X}_n(\cdot)f = f$. In particular, for all $f\in \Cc_c(\Vc^{\rm c})$, $\Ac_{\Gc^{\rm c}}\mathcal{X}_n(\cdot)f\to \Ac_{\Gc^{\rm c}}f$, as $n\to\infty$.

3) Noticing that $[\Delta_{\Gc^{\rm c}},\mathcal{X}_n(\cdot)]=[\Delta_{\Gc_1^{\rm c}},\mathcal{X}_n(\cdot)]\otimes\frac{1}{m_2}$, a straightforward computation ensures that there exists $c$ such that
\[\|\Ac_{\Gc^{\rm c}}[\Delta_{\Gc^{\rm c}},\mathcal{X}_n(\cdot)]\|\leq \frac{c}{\langle n \rangle}.\]
Finally for all $z\in \C\setminus\R$, the condition $\mathcal{X}_n(\cdot)(\Delta_{\Gc^{\rm c}}-z)^{-1}\Cc_c(\Vc^{\rm c})\subset \Cc_c(\Vc^{\rm c})$ is immediate as $\mathcal{X}_n$ is with finite support. \cite[Lemma A.2]{GoMo} gives that $\Delta_{\Gc^{\rm c}}^{\rm le}\in \Cc^1(\Ac_{\Gc^{\rm c}}^{\rm le})$. \qed

\begin{lemma}\label{l:gec}
We have  $e^{\rmi t \Ac_{\Gc^{\rm c}}} \Dc(\Delta_{\Gc^{\rm c}})\subset\Dc(\Delta_{\Gc^{\rm c}})$ for all $t\in \R$.
\end{lemma}
\proof We have $\Delta_{\Gc^{\rm c}}\in \Cc^1(\Ac_{\Gc^{\rm c}})$ and $[\Delta_{\Gc^{\rm c}},\rmi \Ac_{\Gc^{\rm c}}]_\circ$ is bounded. Therefore \cite{GeGe} gives the result. \qed

\begin{lemma}\label{l:C2cusp}
We have $\Delta_{\Gc^{\rm c}}\in \Cc^2(\Ac_{\Gc^{\rm c}})$ and
\begin{align} \label{e:C2cusp}[[\Delta_{\Gc^{\rm c}},\rmi \Ac_{\Gc^{\rm c}}]_\circ,\rmi \Ac_{\Gc^{\rm c}}]_\circ\simeq[[\Delta_{\Gc^{\rm c}}^{\rm le}, \rmi \Ac_{\Gc^{\rm c}}^{\rm le}]_\circ, \rmi \Ac_{\Gc^{\rm c}}^{\rm le}]_\circ\oplus 0.
\end{align}
\end{lemma}
\proof Recalling \eqref{e:C1cusp} and Lemma \ref{Lemma:2}, the result follows from noticing that
$w^{\rm c}(\Delta_{\rm \Gc^{\rm c}}^{\rm le})$ is in $\Cc^1(\Ac_{\Gc^{\rm c}}^{\rm le})$ as product of bounded elements of $\Cc^1(\Ac_{\Gc^{\rm c}}^{\rm le})$. \qed

Concerning the Mourre estimate, we prove the following result:
\begin{proposition}
We have $\Delta_{\Gc^{\rm c}}\in\Cc^2(\Ac_{\Gc^{\rm c}})$. Given a compact interval $\Ic\subset \left(\frac{\alpha}{m_2},\frac{\beta}{m_2} \right)$,
there are $c>0$, a compact operator $K$ such that
 \begin{align}
\label{e:mourrec}
 E_\Ic(\Delta_{\Gc^{\rm c}})[\Delta_{\Gc^{\rm c}}, \rmi \Ac_{\Gc^{\rm c}}]_\circ E_\Ic(\Delta_{\Gc^{\rm c}})&\geq c E_\Ic(\Delta_{\Gc^{\rm c}}) + K,\end{align} in the form sense.
\end{proposition}
\proof
The Lemma \ref{l:C2cusp} provides that $\Delta_{\Gc^{\rm c}}\in \Cc^2(\Ac_{\Gc^{\rm c}})$. On $\Hc^{\rm he}$, $E_\Ic(\Delta^{\rm he}_{\Gc^{\rm c}})$ is compact
since $\Delta^{\rm he}_{\Gc^{\rm c}}$ is with compact resolvent and $\Ic$ is with compact support. With respect to \eqref{e:lhdecomp}, we have $E_\Ic(\Delta_{\Gc^{\rm c}})= E_\Ic(\Delta_{\Gc^{\rm c}}^{\rm le})\oplus E_\Ic(\Delta_{\Gc^{\rm c}}^{\rm he})$ and
 \begin{align*}
 E_\Ic(\Delta_{\Gc^{\rm c}})[\Delta_{\Gc^{\rm c}}, \rmi \Ac_{\Gc^{\rm c}}]_\circ E_\Ic(\Delta_{\Gc^{\rm c}})&=
       E_\Ic(\Delta^{\rm le}_{\Gc^{\rm c}})[\Delta^{\rm le}_{\Gc^{\rm c}}, \rmi \Ac^{\rm le}_{\Gc^{\rm c}}]_\circ E_\Ic(\Delta^{\rm le}_{\Gc^{\rm c}}) \oplus 0
       \\
       &\geq
       cE_\Ic(\Delta^{\rm le}_{\Gc^{\rm c}})\oplus 0\geq c E_\Ic(\Delta_{\Gc^{\rm c}}) + K,
        \end{align*}
in the form sense, where $K$ is a compact operator and
\[c:=\frac{m_2}{2}\inf_{x\in I}\left(x-\frac{\alpha}{m_2}\right)\left(\frac{\beta}{m_2}-x\right)>0.\]
This concludes the proof.\qed

To lighten the text we did not expand more consequences of the Mourre theory in this case and refer to Theorem \ref{t:mainth} for them.

\subsection{The compact part}\label{section5}
We define the conjugate operator on $\ell^2(\Vc)=\ell^2(\Vc^{\rm f})\oplus  \ell^2(\Vc^{\rm 0})\oplus  \ell^2(\Vc^{\rm c})$ as
\begin{align*}
\ell^2(\Vc)&=\ell^2(\Vc^{\rm f})\oplus  \ell^2(\Vc^{\rm 0})\oplus  \ell^2(\Vc^{\rm c})
\\
\Ac&:= \Ac_{\Gc^{\rm f}} \oplus 0 \oplus \Ac_{\Gc^{\rm c}}.
\end{align*}
Since $\Vc^0$ is finite, we have a finite rank perturbation and we conclude that $\Ac$ is self-adjoint and essentially self-adjoint on $\Cc_c(\Vc)$.

\begin{lemma}\label{l:kcC2}
We have $\Delta_\Gc \in \Cc^2(\Ac)$.
\end{lemma}
\proof We have $(\Delta_\Gc - \Delta_{\Gc^{\rm f}}\oplus 0 \oplus \Delta_{\Gc^{\rm c}} )$ that are with finite support. Hence it belongs to $\Cc^2(\Ac)$ by Lemma \ref{Lemma:2}. Next recalling Lemma \ref{l:C2Af} and Lemma \ref{l:C2cusp} we obtain the result. \qed


%




\subsection{The whole graph}\label{section6}
In this section, we give the Mourre estimate in the whole graph.

\begin{proposition}\label{p:mourredelta}
We have $\Delta_\Gc\in \Cc^2(\Ac)$. Given a compact interval $\Ic\subset (\frac{\alpha}{m_2},\frac{\beta}{m_2})$
Moreover, there are $c>0$, a compact operator $K$ such that
\begin{align}\label{e:mourredelta}
E_\Ic(\Delta_\Gc)[\Delta_\Gc, \rmi \Ac]_\circ E_\Ic(\Delta_\Gc)\geq c E_\Ic(\Delta_\Gc) +K.
\end{align}
\end{proposition}
\proof First $\Delta_\Gc\in \Cc^2(\Ac)$ by Lemma \ref{l:kcC2}. Then by collecting \eqref{e:mourrec} and \eqref{e:mourref}, we obtain
\begin{align*}
E_\Ic(\Delta_{\Gc^{\rm f}}\oplus 0 \oplus \Delta_{\Gc^{\rm c}})[\Delta_{\Gc^{\rm f}}\oplus 0 \oplus \Delta_{\Gc^{\rm c}}, \rmi \Ac]_\circ E_\Ic(\Delta_{\Gc^{\rm f}}\oplus 0 \oplus \Delta_{\Gc^{\rm c}})&
\\
&\hspace*{-3cm}\geq c E_\Ic(\Delta_{\Gc^{\rm f}}\oplus 0 \oplus \Delta_{\Gc^{\rm c}}) +K.
\end{align*}
Since the operators $\Delta_\Gc$ and $\Delta_{\Gc^{\rm f}} \oplus 0 \oplus \Delta_{\Gc^{\rm c}}$ are in $\Cc^1_u(\Ac)$ (as in $\Cc^2(A)$, see \cite{ABG}),  \cite[Theorem 7.2.9]{ABG} implies \eqref{e:mourredelta}. \qed

\section{The perturbed model}\label{section7}
In this section, we perturb the metrics of the previous case which will be small to infinity. We obtain similar results however the proof is more involved because we rely on the optimal class  $\Cc^{1,1}(\Ac)$ of the Mourre theory.
\subsection{Perturbation of the metric}\label{Model}
Let $\Gc_{\varepsilon, \mu}:=(\Vc, \Ec_{\varepsilon}, m_\mu)$ where
\[m_\mu(x):=(1+\mu(x))m(x)  \mbox{ and } \Ec_\varepsilon(x,y):=(1+\varepsilon(x,y))\Ec(x,y),\]
where
 $\mu>-1$, $\varepsilon >-1$,
and
\begin{align}\label{e:tend0}
\mu(x)\rightarrow 0  \text{ if } |x|\rightarrow\infty \text{ and } \varepsilon(x,y)\rightarrow0 \text{ if } |x|,|y|\rightarrow\infty.
\end{align}
We set
\[m_\mu^{\ast}:={m_\mu}{|_{\Vc^{\ast}}},\quad  \Ec_\varepsilon^{\ast}:=\Ec_\varepsilon{|_{\Vc^{\ast}\times\Vc^{\ast}}},\]
$\mu^{\ast}:=\mu^{\ast}{|_{\Vc^{\ast}}}$, and $\varepsilon^{\ast}:= \varepsilon_{|_{\Vc^{\ast}\times\Vc^{\ast}}}$, with $\ast\in\{\rm c, \rm f\}$.

To analyse the spectral properties of $\Delta_{\Gc_{\varepsilon, \mu}}$, we compare it to $\Delta_\Gc$. As they do not act in the same spaces, we rely on Proposition \ref{p:uni}.
and send $\Delta_{\Gc_{\varepsilon, \mu}}$ in $\ell^2(\Vc,m)$ with the help of the unitary transformation. Namely, supposing \eqref{e:tend0}. Let
 \begin{align*}
 \widetilde{\Delta}_{\Gc_{\varepsilon, \mu}}&:= T_{m_{\mu}\rightarrow m}\Delta_{\Gc_{\varepsilon, \mu}}T^{-1}_{m_{\mu}\rightarrow m}.
 \end{align*}

A straightforward calculus ensures:
\begin{lemma} 
For all $f\in\Cc_c(\Vc)$, we have
\begin{align}
&\nonumber(\widetilde\Delta_{\Gc_{\varepsilon, \mu}}-\Delta_{\Gc})f(x):=
\frac{1}{m(x)}\sum_{y\sim x}\Bigg(\frac{\varepsilon(x,y)}{\sqrt{(1+\mu(x))(1+\mu(y))}}
\\
\label{e:difr}
&\quad-\frac{\mu(x)+\mu(y)+\mu(x)\mu(y)}{\sqrt{(1+\mu(x))(1+\mu(z))}(1+\sqrt{(1+\mu(x))(1+\mu(z))})}\Bigg)\Ec(x,z)\left(f(x)-f(y)\right)
\\
&\nonumber-\frac{1}{m(x)}\sum_{z\sim x}(1+\varepsilon(x,z))\Ec(x,z)
\frac{\mu(z)-\mu(x)}{(1+\mu(x))\sqrt{1+\mu(z)}(\sqrt{1+\mu(z)}+\sqrt{1+\mu(x)})}f(x).\end{align}
\end{lemma}

\begin{proposition}\label{com}Let $V:\Vc\to\R$ be a function, obeying $V(x)\to0$ if $|x|\to\infty$. We assume that \eqref{e:tend0} holds true, 
then $\widetilde\Delta_{\Gc_{\varepsilon, \mu}}-\Delta_{\Gc}\in \Kc(\ell^2(\Vc),m)$. 
 In particular
\begin{enumerate}
\item[(1)]$\Dc(\Delta_{\Gc_{\varepsilon, \mu}}+V(\cdot))=\Dc(T_{m_{\mu}\rightarrow m}^{-1}\Delta_{\Gc}T_{m_{\mu}\rightarrow m})$,
\item[(2)]$\Delta_{\Gc_{\varepsilon, \mu}}+V(\cdot)$ is essentially self-adjoint on $\Cc_c(\Vc)$,
\item[(3)]$\sigma_{\rm ess}(\Delta_{\Gc_{\varepsilon, \mu}}+V(\cdot))=\sigma_{\rm ess}(\Delta_{\Gc})$.
\end{enumerate}
\end{proposition}
\proof Use Propositions \ref{com.fun} and \ref{com.cusp} and note that the contribution arising from $\Vc^0$ is a finite rank perturbation. \qed

\subsection{Main result}\label{P.result}
 The main result of this section is the following theorem:
\begin{theorem}\label{t:mainth}
Let $\Gc_{\varepsilon, \mu}$ a graph satisfies a condition \eqref{e:tend0} and
\[\Ac_{\Gc_{\varepsilon, \mu}}:=\Ac_{\Gc_{\varepsilon, \mu}^{\rm f}}\oplus 0\oplus\Ac_{\Gc_{\varepsilon, \mu}^{\rm c}}\] 
 be a self-adjoint operator, where $\Ac_{\Gc_{\varepsilon, \mu}^{\ast}}:=T_{m_{\mu^{\ast}}\rightarrow m}^{-1}\Ac_{\Gc^{\ast}}T_{m_{\mu^{\ast}}\rightarrow m}$  
 with $\ast\in \{\rm f,\rm c\}$. Let $V:\Vc\to\R$ be a function such that $V, \varepsilon$, and $\mu$ are radial on $\Vc^{\rm c}$(see Definition \ref{d:rad}).
We assume that:
 \begin{align*}&(H1) \sup_{(x_1,x_2)\in\Vc^{\ast}}\langle x_1\rangle^{1+\epsilon} |V(x_1-1,x_2)-V(x_1,x_2)|<\infty,
 \\
 &(H2) \sup_{(x_1,x_2)\in\Vc^{\ast}} \langle x_1\rangle^{\epsilon+1}|\mu^{\ast}(x_1-1,x_2)-\mu^{\ast}(x_1,x_2)|<\infty,
 \\
 &(H3) \sup_{(x_1,x_2)\in\Vc^{\ast}} \langle x_1\rangle^{\epsilon+1}|\varepsilon^{\ast}((x_1,x_2),(x_1+1,x_2))-\varepsilon^{\ast}((x_1-1,x_2),(x_1,x_2))|<\infty, \end{align*}
 where $V(x)\to0$ if $|x|\to\infty$.
 Then $\Delta_{\Gc_{\varepsilon, \mu}}+V(\cdot)\in \Cc^{1,1}(\Ac_{\Gc_{\varepsilon, \mu}})$. Moreover, for all compact interval $I\subset (\frac{\alpha}{m_2},\frac{\beta}{m_2})$, with $\alpha, \beta$ are given in \eqref{e:c1c2}, there are $c>0$ and a compact operator $K$ such that
\begin{align}\nonumber
&E_I(\Delta_{\Gc_{\varepsilon, \mu}}+V(\cdot))[\Delta_{\Gc_{\varepsilon, \mu}}+V(\cdot), \rmi \Ac_{\Gc_{\varepsilon, \mu}}]_\circ E_I(\Delta_{\Gc_{\varepsilon, \mu}}+V(\cdot))
\\
&\label{Mour}\hspace*{+3.8cm}\geq c E_I(\Delta_{\Gc_{\varepsilon, \mu}}+V(\cdot)) +K.
\end{align} in the form sense. Set $\kappa(\Delta_{\Gc_{\varepsilon, \mu}}+V(\cdot)):=\sigma_{p}(\Delta_{\Gc_{\varepsilon, \mu}}+V(\cdot))\cup\{\frac{\alpha}{m_2},\frac{\beta}{m_2}\}$ where $\sigma_p$ denotes the pure point spectrum. Take $s>1/2$ and $[a,b]\subset \R\setminus \kappa(\Delta_{\Gc_{\varepsilon, \mu}}+V(\cdot))$. We obtain:
\begin{enumerate}
\item[(2)] The eigenvalues of $\Delta_{\Gc_{\varepsilon, \mu}}+V(\cdot)$ distinct from $\alpha$ and $\beta$ are of finite multiplicity and can accumulate
only toward $\alpha$ and $\beta$.
\item[(3)] The singular continuous spectrum of $\Delta_{\Gc_{\varepsilon, \mu}}+V(\cdot)$  is empty.
\item[(4)] The following limit exists and finite:
\[\lim_{\rho\to0}  \sup_{\lambda \in [a,b]}\|\langle \Lambda\rangle^{-s}(\Delta_{\Gc_{\varepsilon, \mu}}+V(\cdot)-\lambda-\rmi\rho)^{-1}\langle \Lambda\rangle^{-s}\|<\infty,\]
\item[(5)] There exists $c>0$ such that for all $f\in\ell^2(\Vc, m_\mu)$,  we have:
\[\int_{\R}\|\langle \Lambda\rangle^{-s}e^{-\rmi t(\Delta_{\Gc_{\varepsilon, \mu}}+V(\cdot))}E_{[a,b]}(\Delta_{\Gc_{\varepsilon, \mu}}+V(\cdot))f\|^2dt\leq c\|f\|^2,\]
with $\Lambda:=\Lambda^{\rm f}\oplus0\oplus\Lambda^{\rm c}$.
\end{enumerate}
\end{theorem}
 \proof First $\Delta_{\Gc_{\varepsilon, \mu}}+V(\cdot)\in \Cc^{1,1}(\Ac_{\Gc_{\varepsilon, \mu}})$ because $\Delta_{\Gc_{\varepsilon, \mu}^{\rm f}}\oplus 0 \oplus \Delta_{\Gc_{\varepsilon, \mu}^{\rm c}}\in \Cc^{1,1}(\Ac_{\Gc_{\varepsilon, \mu}})$ by the Lemma \ref{C1,1fun}, the Lemma \ref{C1,1cusp}, 
 the Lemma  \ref{l:Vcusp}, the Lemma \ref{l:Vfun} and by Lemma \ref{l:kcC2}. In particular, we have that the two operators are in $\Cc^{1}_u(\Ac_{\Gc_{\varepsilon, \mu}})$, see \cite{ABG}.

Then, using the Proposition \ref{mourfun2} and the Proposition \ref{mourcusp2} 
we obtain
\begin{align*}
E_I(\Delta_{\Gc_{\varepsilon, \mu}^{\rm f}}\oplus 0 \oplus \Delta_{\Gc_{\varepsilon, \mu}^{\rm c}}+V(\cdot))[\Delta_{\Gc_{\varepsilon, \mu}^{\rm f}}\oplus 0 \oplus &\Delta_{\Gc_{\varepsilon, \mu}^{\rm c}}+V(\cdot), \rmi \Ac_{\Gc_{\varepsilon, \mu}}]_\circ
\\
&\hspace*{1cm}E_I(\Delta_{\Gc_{\varepsilon, \mu}^{\rm f}}\oplus 0 \oplus \Delta_{\Gc_{\varepsilon, \mu}^{\rm c}}+V(\cdot))
\\
&\geq c E_I(\Delta_{\Gc_{\varepsilon, \mu}^{\rm f}}\oplus 0 \oplus \Delta_{\Gc_{\varepsilon, \mu}^{\rm c}}+V(\cdot)) +K.
\end{align*}
Since $\Delta_{\Gc_{\varepsilon, \mu}^{\rm f}}\oplus 0 \oplus \Delta_{\Gc_{\varepsilon, \mu}^{\rm c}} \in \Cc^1_u(\Ac_{\Gc_{\varepsilon, \mu}})$, $\Delta_{\Gc_{\varepsilon, \mu}^{\rm f}}\oplus 0 \oplus \Delta_{\Gc_{\varepsilon, \mu}^{\rm c}}-\Delta_{\Gc_{\varepsilon, \mu}}\in \Kc(\ell^2(\Vc,m_\mu))$, 
and $V(\cdot)\in\Cc^1_u(\Ac_{\Gc_{\varepsilon, \mu}})$ and by \cite[Theorem 7.2.9]{ABG}, we obtain \eqref{Mour}.
By Lemma \ref{l:Vfun} and Lemma \ref{l:Vcusp}, $V(\cdot)\in\Cc^{1,1}(\Ac_{\Gc_{\varepsilon, \mu}})$. And by using Proposition \ref{com.fun} and Proposition \ref{com.cusp}, we have that $(\Delta_{\Gc_{\varepsilon, \mu}^{\rm f}}\oplus 0 \oplus \Delta_{\Gc_{\varepsilon, \mu}^{\rm c}}+\rm i)^{-1}-(\Delta_{\Gc_{\varepsilon, \mu}}+\rm i)^{-1}\in\Kc(\ell^2(\Vc,m_\mu))$. 
Finally, we turn to points (4). It is enough to obtain them with $s\in (1/2, 1)$.
We apply \cite[Proposition 7.5.6]{ABG} and obtain
\[\lim_{\rho\to0}\|\langle \Ac\rangle^{-s}(\Delta_{\Gc_{\varepsilon, \mu}}+V(\cdot)-\lambda-\rm i\rho)^{-1}\langle \Ac\rangle^{-s}\|, \]
exists and finite.
Using Propositions \ref{p:hypothesis} b) and \ref{c:hypothesis} b)
\[\|\langle \Ac_{\Gc_{\varepsilon, \mu}}\rangle f\|\leq a\|\langle \Lambda\rangle f\|, \]
for all $f\in\Dc(\Lambda)$. By Riesz-Thorin interpolation, there is $a_s>0$ such that
\[\|\langle \Ac_{\Gc_{\varepsilon, \mu}}\rangle^{ s} f\|\leq a_s\|\langle \Lambda\rangle^{ s}f\|, \]
for all $f\in\Dc(\Lambda^s)$. We conclude that
 $\lim_{\rho\to0}\|\langle \Lambda\rangle^{-s}(\Delta_{\Gc_{\varepsilon, \mu}}+V(\cdot)-\lambda-\rm i\rho)^{-1}\langle \Lambda\rangle^{-s}\|$  exists and finite. The point (5) is an immediate consequence of (4).  \qed

\subsection{The funnel side}
We first deal with the question of the essential spectrum.
\begin{proposition}\label{com.fun}
Let $V^{\rm f}:\Vc^{\rm f}\to \R$ be a function obeying $V^{\rm f}(x)\to0$ if $|x|\to\infty$. We assume that \eqref{e:tend0} holds true 
then $\widetilde{\Delta}_{\Gc_{\varepsilon, \mu}^{\rm f}}-\Delta_{\Gc^{\rm f}}\in \Kc(\ell^2(\Vc^{\rm f}),m^{\rm f})$, where $\widetilde\Delta_{\Gc_{\varepsilon, \mu}}:= T_{m_{\mu}\rightarrow m}\Delta_{\Gc_{\varepsilon, \mu}}T^{-1}_{m_{\mu}\rightarrow m}$. 
In particular,
\begin{enumerate}
\item[(1)]$\Dc(\Delta_{\Gc_{\varepsilon, \mu}^{\rm f}}+V(\cdot))=\Dc(T_{m_{\mu}\rightarrow m}^{-1}\Delta_{\Gc}T_{m_{\mu}\rightarrow m})$,
\item[(2)]$\Delta_{\Gc_{\varepsilon, \mu}^{\rm f}}+V(\cdot)$ is essentially self-adjoint on $\Cc_c(\Vc)$,
\item[(3)]$\sigma_{\rm ess}(\Delta_{\Gc_{\varepsilon, \mu}^{\rm f}}+V(\cdot))=\sigma_{\rm ess}(\Delta_{\Gc})$.
\end{enumerate}
\end{proposition}

\proof 
We shall show that $\widetilde{\Delta}_{\Gc_{\varepsilon, \mu}^{\rm f}}-\Delta_{\Gc^{\rm f}}\in\Kc(\ell^2(\Vc^{\rm f},m^{\rm f}))$, as in \eqref{e:difr}. Let $f\in \Cc_c(\Vc)$,
\begin{align*}
&\big|\langle f, (\widetilde{\Delta}_{\Gc_{\varepsilon, \mu}^{\rm f}}-\Delta_{\Gc^{\rm f}}) f\rangle_{\ell^2(\Vc^{\rm f},m^{\rm f})}\big|=\bigg|\sum_{x\in\Vc^{\rm f}} m^{\rm f}(x)\left((\widetilde{\Delta}_{\Gc_{\varepsilon, \mu}^{\rm f}}-\Delta_{\Gc^{\rm f}}) f\right)(x)\overline{f(x)}\bigg|
\\
&\leq\bigg|\sum_{x\in\Vc^{\rm f}} m^{\rm f}(x)\frac{1}{m^{\rm f}(x)}\sum_{z\sim x}\frac{\varepsilon^{\rm f}(x,z)}{\sqrt{(1+\mu^{\rm f}(x))(1+\mu^{\rm f}(z))}}\Ec^{\rm f}(x,z)
\\
&\hspace*{+1cm}\times\big(f(x)-f(z)\big)\overline{f(x)}\bigg|
\\
&\quad+\bigg|\sum_{x\in\Vc^{\rm f}} m^{\rm f}(x)\frac{1}{m^{\rm f}(x)}\sum_{z\sim x}\frac{1-\sqrt{(1+\mu^{\rm f}(x))(1+\mu^{\rm f}(z))}}{\sqrt{(1+\mu^{\rm f}(x))(1+\mu^{\rm f}(z))}}\Ec^{\rm f}(x,z)
\\
&\hspace*{+1cm}\times\big(f(x)-f(z)\big)\overline{f(x)}\bigg|
+\big|\langle f,W^{\rm f}(\cdot)f\rangle\big|
\\
&\leq 2\langle f, (\deg_1(\cdot)+\deg_2(\cdot)+|W^{\rm f}(\cdot)|) f\rangle, \end{align*}  with
 \[\deg_1(x):=\frac{1}{m^{\rm f}(x)}\sum_{z\in\Vc^{\rm f}}\frac{\varepsilon^{\rm f}(x,z)}{\sqrt{(1+\mu^{\rm f}(x))(1+\mu^{\rm f}(z))}}\Ec^{\rm f}(x,z)\]
and
 \[\deg_2(x):=\frac{1}{m^{\rm f}(x)}\sum_{z\in\Vc^{\rm f}}\left|\frac{1-\sqrt{(1+\mu^{\rm f}(x))(1+\mu^{\rm f}(z))}}{\sqrt{(1+\mu^{\rm f}(x))(1+\mu^{\rm f}(z))}}\right|\Ec^{\rm f}(x,z), \]
for all $x=(x_1,x_2)\in\Vc^{\rm f}$.
We have \begin{align*}|\deg_1(x)|&=\left|\frac{1}{m^{\rm f}(x)}\sum_{z\in\Vc^{\rm f}}\frac{\varepsilon^{\rm f}(x,z)}{\sqrt{(1+\mu^{\rm f}(x))(1+\mu^{\rm f}(z))}}\Ec^{\rm f}(x,z)\right|
\\
&\leq\sup_{z\sim x}\left|\frac{\varepsilon^{\rm f}(x,z)}{\sqrt{(1+\mu^{\rm f}(x))(1+\mu^{\rm f}(z))}}\right|\deg_{\Gc^{\rm f}}(x).
\end{align*}
Since $\Vc_2$ is a finite set and for all $x_2\in\Vc_2$, $\frac{\varepsilon^{\rm f}((x_1,x_2),(z_1,z_2))}{\sqrt{(1+\mu^{\rm f}(x_1,x_2))(1+\mu^{\rm f}(z_1,z_2))}}\to 0$  when $x_1,z_1\to\infty$ and since $\deg_{\Gc^{\rm f}}(\cdot)$ is bounded then $\deg_1(\cdot)$ is compact. In the same way, using that $\forall x_2,z_2\in\Vc_2,~\frac{1-\sqrt{(1+\mu^{\rm f}(x_1,x_2))(1+\mu^{\rm f}(z_1,z_2))}}{\sqrt{(1+\mu^{\rm f}(x_1,x_2))(1+\mu^{\rm f}(z_1,z_2))}}\to 0$ if $x_1,z_1\to\infty$, we obtain the compactness of $\deg_2(\cdot)$. 

Now, we will show that $W^{\rm f}\in\Kc(\ell^2(\Vc^{\rm f},m^{\rm f}))$. For all  $x\in\Vc^{\rm f},$  we have \begin{align*}|W^{\rm f}(x)|&=\bigg|\frac{1}{m^{\rm f}(x)}\sum_{z\sim x}(1+\varepsilon^{\rm f}(x,z))\Ec^{\rm f}(x,z)
\\
&\hspace*{-0.2cm}\quad\times\left(\frac{\mu^{\rm f}(z)-\mu^{\rm f}(x)}{(1+\mu^{\rm f}(x))\sqrt{1+\mu^{\rm f}(z)}(\sqrt{1+\mu^{\rm f}(z)}+\sqrt{1+\mu^{\rm f}(x)})}\right)\bigg|
\\
& \hspace*{-0.2cm}
\leq\sup_{z\sim x}\left|(1+\varepsilon^{\rm f}(x,z))\left(\frac{\mu^{\rm f}(z)-\mu^{\rm f}(x)}{(1+\mu^{\rm f}(x))\sqrt{1+\mu^{\rm f}(z)}(\sqrt{1+\mu^{\rm f}(z)}+\sqrt{1+\mu^{\rm f}(x)})}\right)\right|
\\
& \hspace*{-0.2cm} \quad\times\deg_{\Gc^{\rm f}}(x).\end{align*} Since $\Vc_2$ is a finite set and $\big(1+\varepsilon^{\rm f}(x,z)\big)\big(\mu^{\rm f}(z)-\mu^{\rm f}(x)\big)\to0$ when $|x|,|z|\to\infty$, $\deg_{\Gc^{\rm f}}(\cdot)$ is bounded  and since $V^{\rm f}(\cdot)$ is a compact perturbation, we conclude that $\widetilde{\Delta}_{\Gc_{\varepsilon, \mu}^{\rm f}}-\Delta_{\Gc^{\rm f}}$ is compact. The points (1) and (2) follow from Theorem \cite[Theorem XIII.14]{RS} and (3) from the Weyl's Theorem. \qed

We turn to the Mourre estimate.

\begin{proposition}\label{mourfun2} Let $V^{\rm f}:\Vc^{\rm f}\to\R$ be a function. We assume that $(H1)$, $(H2)$, and $(H3)$ hold true,
 where $\varepsilon^{\rm f}(x,z)\to0$ if $|x|,|z|\to\infty$, $\mu^{\rm f}(x)\to0$ if $|x|\to\infty$ and $V^{\rm f}(x)\to0$ if $|x|\to\infty$.
 Then $\Delta_{\Gc_{\varepsilon, \mu}^{\rm f}}+V^{\rm f}(\cdot)\in\Cc^{1.1}(\Ac_{\Gc_{\varepsilon, \mu}^{\rm f}})$. Moreover, for all compact interval $I\subset (\frac{\alpha}{m_2},\frac{\beta}{m_2})$, there are $c>0$, a compact operator $K$ such that
 \begin{align}
 \nonumber
 E_I(\Delta_{\Gc_{\varepsilon, \mu}^{\rm f}}+V^{\rm f}(\cdot))[\Delta_{\Gc_{\varepsilon, \mu}^{\rm f}}+V^{\rm f}(\cdot), \rmi \Ac_{\Gc_{\varepsilon, \mu}^{\rm f}}]_\circ E_I(\Delta_{\Gc_{\varepsilon, \mu}^{\rm f}}+V^{\rm f}(\cdot))&
 \\
 \label{mourrF}
 &\hspace*{-3cm} \geq c E_I(\Delta_{\Gc_{\varepsilon, \mu}^{\rm f}}+V^{\rm f}(\cdot)) + K,\end{align} in the form sense.
\end{proposition}
\proof The Proposition \ref{C1,1fun} and Lemma \ref{l:Vfun} give that ${\Delta}_{\Gc_{\varepsilon, \mu}^{\rm f}}+V^{\rm f}(\cdot)\in\Cc^{1,1}(\Ac_{\Gc_{\varepsilon, \mu}^{\rm f}})$. Since $\widetilde{\Delta}_{\Gc_{\varepsilon, \mu}^{\rm f}}-\Delta_{\Gc^{\rm f}}$ is a compact operator by Proposition \ref{com.fun}, thanks to \eqref{e:mourref} and by \cite[Theorem 7.2.9]{ABG}, we obtain \eqref{mourrF}.\qed

We start with a technical lemma so as to apply \cite[Proposition 7.5.7]{ABG}.

\begin{proposition}\label{p:hypothesis}
Let $\Lambda^{\rm f}:=(Q+1/2)\otimes\un_{\Vc^{\rm f}}$. It satisfies the following assertions:
\begin{enumerate}
\item[(1)] $e^{\rmi\Lambda^{\rm f}t }\Dc(\Delta_{\Gc_{\varepsilon, \mu}^{\rm f}})\subset \Dc(\Delta_{\Gc_{\varepsilon, \mu}^{\rm f}})$  and there exists a finite constant $\rm c$, such that  \[\|e^{\rmi\Lambda^{\rm f}t}\|_{\Bc(\Dc(\Delta_{\Gc_{\varepsilon, \mu}^{\rm f}}))} \leq \rm c, \quad \hbox{for all}\quad t\in\R.\]
\item[(2)] $\Dc(\Lambda^{\rm f})\subset \Dc(\Ac_{\Gc_{\varepsilon, \mu}^{\rm f}})$.
\item[(3)] $(\Lambda^{\rm f})^{-2}(\Ac_{\Gc_{\varepsilon, \mu}^{\rm f}})^2$ extends to a continuous operator in $\Dc(\Delta_{\Gc_{\varepsilon, \mu}^{\rm f}})$.
\end{enumerate} Note that $\Delta_{\Gc_{\varepsilon, \mu}^{\rm f}}$ is bounded then $\Dc(\Delta_{\Gc_{\varepsilon, \mu}^{\rm f}})=\ell^2(\Vc^{\rm f}, m_\mu^{\rm f})$.
\end{proposition}
\proof With the help of the unitary transformation $T_{m_\mu\to m}$, it is enough to prove the result with $\varepsilon=0$ and $\mu=0$.

(1) Since $\Delta_{\Gc^{\rm f}}$ is bounded it is verified by a functional calculus.

(2) Let $f\in\Cc_c(\Vc^{\rm f})$ ,
\begin{align*}
&\|\Ac_{\Gc^{\rm f}}f\|^2_{\ell^2(\Vc^{\rm f},m^{\rm f})}
\\
&=\sum_{x\in\Vc^{\rm f}}m^{\rm f}(x)\left|\frac{\rmi}{2}\left(e^{1/2}(Q-1/2)U\otimes\un_{\Vc_2^{\rm f}}-e^{-1/2}(Q+1/2)U^*\otimes\un_{\Vc_2^{\rm f}}\right)f(x)\right|^2
\\
&\leq c \sum_{x\in\Vc^{\rm f}}m^{\rm f}(x)\left|(Q+1/2)\otimes\un_{\Vc_2^{\rm f}}f(x)\right|^2 \leq c\|\Lambda^{\rm f}f\|^2_{\ell^2(\Vc^{\rm f},m^{\rm f})}.
\end{align*}
Since $\Lambda^{\rm f}$ is essentially self-adjoint, the result follows.

(3) For all $f\in\Cc_c(\Vc^{\rm f})$, and by using the relations of Subsection \ref{dec}, we have
\begin{align*}&\Ac^2_{m_1}f(n)=\frac{1}{4}(2n^2+1/2)f(n)-\frac{1}{4}e(n-1/2)(n-3/2)f(n-2)
\\
&\quad-\frac{1}{4}e^{-1}(n+1/2)(n+3/2)f(n+2).\end{align*}
Then for all $f\in\Cc_c(\Vc^{\rm f})$.
\begin{align*}&\|(\Lambda^{\rm f})^{-2}(\Ac_{\Gc^{\rm f}})^2f\|^2
=
\\
&\sum_{(x_1,x_2)\in\Vc^{\rm f}}m^{\rm f}(x_1,x_2)\bigg|\frac{1}{4}\left(\Big((Q+1/2)^{-2}(2Q^2+1/2)\Big)\otimes\un_{\Vc_2^{\rm f}}\right)f(x_1,x_2)
\\
&\quad-\frac{1}{4}e\left(\Big((Q+1/2)^{-2}(Q-1/2)(Q-3/2)\Big)\otimes\un_{\Vc_2^{\rm f}}\right)f(x_1-2,x_2)
\\
&\quad-\frac{1}{4}e^{-1}\left(\Big((Q+1/2)^{-2}(Q+1/2)(Q+3/2)\Big)\otimes\un_{\Vc_2^{\rm f}}\right)f(x_1+2,x_2)\bigg|^2.
\end{align*}
Then, there exists $C>0$ such that for all $f\in\Cc_c(\Vc^{\rm f})$, $\|(\Lambda^{\rm f})^{-2}(\Ac_{\Gc^{\rm f}})^2f\|^2\leq C \|f\|^2$. By density, we find the result. \qed

The proof of Proposition \ref{C1,1fun} will be long and technical. For the sake of the reader, we have separated the treatment of the potential $V^{\rm f}$ to present the technical steps.

\begin{lemma}\label{l:Vfun}Let $V^{\rm f}:\Vc^{\rm f}\to\R$ be a function. We assume that $(H1)$ holds true, 
then $V^{\rm f}(\cdot)\in \Cc^1(\Ac_{\Gc_{\varepsilon, \mu}^{\rm f}})$ and
 $[V^{\rm f}(\cdot),\Ac_{\Gc_{\varepsilon, \mu}^{\rm f}}]_\circ\in\Cc^{0,1}(\Ac_{\Gc_{\varepsilon, \mu}^{\rm f}})$. In particular, $V^{\rm f}(\cdot)\in\Cc^{1,1}(\Ac_{\Gc_{\varepsilon, \mu}^{\rm f}})$.
\end{lemma}
\proof 
First, recalling $\left[V^{\rm f}(\cdot),\rmi \Ac_{\Gc_{\varepsilon, \mu}^{\rm f}}\right]_\circ=T_{m_{\mu^{\ast}}\rightarrow m}^{-1}\left[V^{\rm f}(\cdot),\rmi \Ac_{\Gc^{\rm f}}\right]_\circ T_{m_{\mu^{\ast}}\rightarrow m}$, it is enough to deal with $\epsilon=\mu=0$.
Next, we recall that
\begin{align*}\left[V^{\rm f}(\cdot),\rmi \Ac_{\Gc^{\rm f}}\right]_\circ&=\frac{e^{-1/2}}{2} \left(Q+\frac{1}{2}\right)\left[V^{\rm f},U^*\right]\otimes\un_{\Vc_2^{\rm f}} +\frac{e^{1/2}}{2} \left(\frac{1}{2}-Q\right)\left[V^{\rm f},U\right]\otimes\un_{\Vc_2^{\rm f}}.
\end{align*}
 By using (H1) at the last step, there is $C$ such that, for all $f\in\Cc_c(\Vc)$,   \begin{align*}
 \big\|\langle Q+1/2\rangle^{\epsilon}\otimes\un_{\Vc^{\rm f}_2}[V^{\rm f}(\cdot),\rmi \Ac_{\Gc^{\rm f}}]_\circ f\big\|&
 \\
&\hspace*{-5cm}\leq\frac{e^{-1/2}}{2}\big\|\big(\langle Q+1/2\rangle^{\epsilon}(Q+1/2)\left[V^{\rm f},U^*\right]\big)\otimes\un_{\Vc_2^{\rm f}}f \big\|
\\
&\hspace*{-5cm}\quad+\frac{e^{1/2}}{2}\big\|\big(\langle Q+1/2\rangle^{\epsilon} (Q+1/2)\left[V^{\rm f},U\right]\big)\otimes\un_{\Vc_2^{\rm f}}f \big\|
\\
&\hspace*{-5cm}\leq\frac{e^{-1/2}}{2}\big\|\langle \Lambda\rangle^{\epsilon+1}\big(\left[V^{\rm f},U^*\right]\otimes\un_{\Vc_2^{\rm f}}\big)f \big\|+\frac{e^{1/2}}{2}\big\|\langle \Lambda\rangle^{\epsilon+1}\big(\left[V^{\rm f},U\right]\otimes\un_{\Vc_2^{\rm f}}\big)f \big\|
\leq C\|f\|.\end{align*}
Finally thanks to Proposition \ref{p:hypothesis}, we can apply \cite[Proposition 7.5.7]{ABG} and the result follows.\qed

We conclude this section with the most technical part.

\begin{proposition}\label{C1,1fun}
Assuming  $(H2)$ and $(H3)$ hold true, 
we have ${\Delta}_{\Gc_{\varepsilon, \mu}^{\rm f}}\in \Cc^1(\Ac_{\Gc_{\varepsilon, \mu}^{\rm f}})$. Moreover
$[{\Delta}_{\Gc_{\varepsilon, \mu}^{\rm f}},\Ac_{\Gc_{\varepsilon, \mu}^{\rm f}}]_\circ\in\Cc^{0,1}(\Ac_{\Gc_{\varepsilon, \mu}^{\rm f}})$. In particular, ${\Delta}_{\Gc_{\varepsilon, \mu}^{\rm f}}\in\Cc^{1,1}(\Ac_{\Gc_{\varepsilon, \mu}^{\rm f}})$.
\end{proposition}

\proof 
We work in $\ell^2(\Vc^{\rm f}, m^{f})$. First, using the computation below with $\epsilon=0$ and recalling that $\Ac_{\Gc^{\rm f}}\Cc_c(\Vc^{\rm f}) \subset \Cc_c(\Vc^{\rm f})$, we get there is $c>0$ such that
\[\|[\widetilde{\Delta}_{\Gc_{\varepsilon, \mu}^{\rm f}},\Ac_{\Gc^{\rm f}}]f\|_{\ell^2(\Vc^{\rm f},m^{\rm f})} \leq c
\|f\|_{\ell^2(\Vc^{\rm f},m^{\rm f})},\]
for all $f\in \Cc_c(\Vc^{\rm f})$.  By density, we obtain that $\widetilde{\Delta}_{\Gc_{\varepsilon, \mu}^{\rm f}}\in \Cc^1(\Ac_{\Gc^{\rm f}})$.
Next take $\epsilon >0$ and $f\in \Cc_c(\Vc^{\rm f})$. We aim at proving that $[\widetilde{\Delta}_{\Gc_{\varepsilon, \mu}^{\rm f}},\Ac_{\Gc^{\rm f}}]_\circ$ is $\Cc^{0,1}(\Ac_{\Gc^{\rm f}})$.
\begin{align}
\nonumber
&\|\langle\Lambda^{\rm f}\rangle^{\epsilon}[\widetilde{\Delta}_{\Gc_{\varepsilon, \mu}^{\rm f}},\Ac_{\Gc^{\rm f}}]f\|_{\ell^2(\Vc^{\rm f},m^{\rm f})}=\sum_{x\in\Vc^{\rm f}}m^{\rm f}(x)\big|\langle\Lambda^{\rm f}\rangle^{\epsilon}[\widetilde{\Delta}_{\Gc_{\varepsilon, \mu}^{\rm f}},\Ac_{\Gc^{\rm f}}]f(x)\big|^2
\\
\nonumber
&\leq\sum_{x\in\Vc^{\rm f}}m^{\rm f}(x)\bigg|\frac{\rmi}{2}\Big(e^{1/2}\langle Q-1/2\rangle^{1+\epsilon}\otimes\un_{\Vc_2^{\rm f}}\Big)
\\
\label{e:com}
&\hspace*{+.2cm}\times\bigg(\frac{1}{m^{\rm f}(x)}\sum_{z\sim x}\Ec^{\rm f}(x,z)\frac{1+\varepsilon^{\rm f}(x,z)}{\sqrt{(1+\mu^{\rm f}(x))(1+\mu^{\rm f}(z))}}f(x_1-1,x_2)
\\
\nonumber
&\hspace*{+0.2cm}-\frac{1}{m^{\rm f}(x_1-1,x_2)}\sum_{z\sim x}\Ec^{\rm f}((x_1-1,x_2),z)f(z_1-1,z_2)\times\frac{1+\varepsilon^{\rm f}((x_1-1,x_2),z)}{\sqrt{(1+\mu^{\rm f}(x_1-1,x_2)(1+\mu^{\rm f}(z))}}\bigg)\bigg|^2
\\
\nonumber
&\quad+\sum_{x\in\Vc^{\rm f}}m^{\rm f}(x)\bigg|\frac{\rmi}{2}\Big(e^{1/2}\langle Q-1/2\rangle^{1+\epsilon}\otimes\un_{\Vc_2^{\rm f}}\Big)
\\
\nonumber
&\hspace*{+0.2cm}\times\bigg(\frac{1}{m^{\rm f}(x_1-1,x_2)}\sum_{z\sim x}\Ec^{\rm f}((x_1-1,x_2),z)\frac{1+\varepsilon^{\rm f}((x_1-1,x_2),z)}{\sqrt{(1+\mu^{\rm f}(x_1-1,x_2)(1+\mu^{\rm f}(z))}}f(z)
\\
\nonumber
&\hspace*{+0.2cm}-\frac{1}{m^{\rm f}(x)}\sum_{z\sim x}\Ec^{\rm f}(x,z)\frac{1+\varepsilon^{\rm f}(x,z)}{\sqrt{(1+\mu^{\rm f}(x))(1+\mu^{\rm f}(z))}}f(z_1-1,z_2)\bigg)\bigg|^2
\\
\nonumber
&\quad+\sum_{x\in\Vc^{\rm f}}m^{\rm f}(x)\bigg|\frac{\rmi}{2}\Big(e^{-1/2}\langle Q+1/2\rangle^{1+\epsilon}\otimes\un_{\Vc_2^{\rm f}}\Big)
\\
\nonumber
&\hspace*{+0.2cm}\times\bigg(\frac{1}{m^{\rm f}(x_1+1,x_2)}\sum_{z\sim x}\Ec^{\rm f}((x_1+1,x_2),z)\frac{1+\varepsilon^{\rm f}((x_1+1,x_2),z)}{\sqrt{(1+\mu^{\rm f}(x_1+1,x_2)(1+\mu^{\rm f}(z))}}
\\
\nonumber
&\hspace*{+0.2cm}-\frac{1}{m^{\rm f}(x)}\sum_{z\sim x}\Ec^{\rm f}(x,z)\frac{1+\varepsilon^{\rm f}(x,z)}{\sqrt{(1+\mu^{\rm f}(x))(1+\mu^{\rm f}(z))}}\bigg)f(x_1+1,x_2)\bigg|^2
\\
\nonumber
&+\sum_{x\in\Vc^{\rm f}}m^{\rm f}(x)\bigg|\frac{\rmi}{2}\Big(e^{-1/2}\langle Q+1/2\rangle^{1+\epsilon}\otimes\un_{\Vc_2^{\rm f}}\Big)
\\
\nonumber
&\hspace*{+0.2cm}\times\bigg(\frac{1}{m^{\rm f}(x)}\sum_{z\sim x}\Ec^{\rm f}(x,z)\frac{1+\varepsilon^{\rm f}(x,z)}{\sqrt{(1+\mu^{\rm f}(x))(1+\mu^{\rm f}(z))}}f(z_1+1,z_2)
\\
\nonumber
&\hspace*{+0.2cm}-\frac{1}{m^{\rm f}(x_1+1,x_2)}\sum_{z\sim x}\Ec^{\rm f}((x_1+1,x_2),z)\frac{1+\varepsilon^{\rm f}((x_1+1,x_2),z)}{\sqrt{(1+\mu^{\rm f}(x_1+1,x_2)(1+\mu^{\rm f}(z))}}f(z_1,z_2)\bigg)\bigg|^2
\\
\nonumber&+\|\langle\Lambda^{\rm f}\rangle^{\epsilon}[W^{\rm f}(\cdot),\Ac_{\Gc^{\rm f}}]f\|_{\ell^2(\Vc^{\rm f},m^{\rm f})},\end{align}
with
\begin{align}\nonumber
&[W^{\rm f}(\cdot),\Ac_{\Gc^{\rm f}}]f(x)=\frac{\rmi}{2}e^{1/2}\Big((Q-1/2)\otimes\un_{\Vc_2^{\rm f}}\Big)\frac{1}{m^{\rm f}(x)}\sum_{z\sim x}(1+\varepsilon^{\rm f}(x,z))\Ec^{\rm f}(x,z)Uf(x)
\\
&\label{commut}\hspace*{+4cm}\times\left(\frac{\mu^{\rm f}(z)-\mu^{\rm f}(x)}{(1+\mu^{\rm f}(x))\sqrt{1+\mu^{\rm f}(z)}(\sqrt{1+\mu^{\rm f}(z)}+\sqrt{1+\mu^{\rm f}(x)})}\right)\end{align}

\begin{align}
\nonumber
&\hspace*{-0.5cm}-\frac{\rmi}{2}e^{1/2}\Big((Q-1/2)\otimes\un_{\Vc_2^{\rm f}}\Big)\frac{1}{m^{\rm f}(x_1-1,x_2)}\sum_{z\sim x}(1+\varepsilon^{\rm f}((x_1-1,x_2),(z_1-1,z_2))
\\
\nonumber
&\hspace*{-1cm}\times\left(\frac{\mu^{\rm f}(z_1-1,z_2)-\mu^{\rm f}(x_1-1,x_2)}{(1+\mu^{\rm f}(x_1-1,x_2))\sqrt{1+\mu^{\rm f}(z_1-1,z_2)}(\sqrt{1+\mu^{\rm f}(z_1-1,z_2)}+\sqrt{1+\mu^{\rm f}(x_1-1,x_2)})}\right)
\\
\nonumber
&\hspace*{+4cm}\times\Ec^{\rm f}((x_1-1,x_2),(z_1-1,z_2))Uf(x)
\\
\nonumber
&\hspace*{-0.5cm}+\frac{\rmi}{2}e^{-1/2}\Big((Q+1/2)\otimes\un_{\Vc_2^{\rm f}}\Big)\frac{1}{m^{\rm f}(x_1+1,x_2)}\sum_{z\sim x}(1+\varepsilon^{\rm f}((x_1+1,x_2),(z_1+1,z_2))
\\
\nonumber
&\hspace*{-1cm}\times\left(\frac{\mu^{\rm f}(z_1+1,z_2)-\mu^{\rm f}(x_1+1,x_2)}{(1+\mu^{\rm f}(x_1+1,x_2))\sqrt{1+\mu^{\rm f}(z_1+1,z_2)}(\sqrt{1+\mu^{\rm f}(z_1+1,z_2)}+\sqrt{1+\mu^{\rm f}(x_1+1,x_2)})}\right)
\\
\nonumber
&\hspace*{+4cm}\times\Ec^{\rm f}((x_1+1,x_2),(z_1+1,z_2))U^*f(x)
\\
\nonumber
&\hspace*{-0.5cm}-\frac{\rmi}{2}e^{-1/2}\Big((Q+1/2)\otimes\un_{\Vc_2^{\rm f}}\Big)\frac{1}{m^{\rm f}(x)}\sum_{z\sim x}(1+\varepsilon^{\rm f}(x,z))\Ec^{\rm f}(x,z)U^*f(x)
\\
\nonumber
&\hspace*{+4cm}\times\left(\frac{\mu^{\rm f}(z)-\mu^{\rm f}(x)}{(1+\mu^{\rm f}(x))\sqrt{1+\mu^{\rm f}(z)}(\sqrt{1+\mu^{\rm f}(z)}+\sqrt{1+\mu^{\rm f}(x)})}\right).\end{align}
 We treat the first term of $\|\langle\Lambda^{\rm f}\rangle^{\epsilon}[\widetilde{\Delta}_{\Gc_{\varepsilon, \mu}^{\rm f}},\Ac_{\Gc^{\rm f}}]f\|_{\ell^2(\Vc^{\rm f},m^{\rm f})}$ in \eqref{e:com}.
\begin{align}\nonumber&\sum_{x\in\Vc^{\rm f}}m^{\rm f}(x)\bigg|\frac{\rmi}{2}\Big(e^{1/2}\langle Q-1/2\rangle^{1+\epsilon}\otimes\un_{\Vc_2^{\rm f}}\Big)\sum_{z\sim x}\bigg(\frac{\Ec^{\rm f}(x,z)(1+\varepsilon^{\rm f}(x,z))}{ m^{\rm f}(x)\sqrt{(1+\mu^{\rm f}(x))(1+\mu^{\rm f}(z))}}
\\
\nonumber
&\hspace*{+0.5cm}-\frac{\Ec^{\rm f}((x_1-1,x_2),z)(1+\varepsilon^{\rm f}((x_1-1,x_2),z))}{m(x_1-1,x_2)\sqrt{(1+\mu^{\rm f}(x_1-1,x_2)(1+\mu^{\rm f}(z))}}\bigg)f(x_1-1,x_2)\bigg|^2\end{align}
\begin{align}
\nonumber
&\leq2\sum_{x\in\Vc^{\rm f}}m^{\rm f}(x)\bigg|\frac{\rmi}{2}\Big(e^{1/2}\langle Q-1/2\rangle^{1+\epsilon}\otimes\un_{\Vc_2^{\rm f}}\Big)
\\
\label{e:part1}
&\hspace*{+0.5cm}\times\sum_{z_1\sim x_1}\delta_{z_2=x_2}\bigg(\frac{\Ec^{\rm f}_1(x_1,z_1)(1+\varepsilon^{\rm f}(x,z))}{ m(x)\sqrt{(1+\mu^{\rm f}(x))(1+\mu^{\rm f}(z))}}
\\
\nonumber
&\hspace*{+0.5cm}-\frac{\Ec^{\rm f}_1(x_1-1,z_1)(1+\varepsilon^{\rm f}((x_1-1,x_2),z))}{m^{\rm f}(x_1-1,x_2)\sqrt{(1+\mu^{\rm f}(x_1-1,x_2)(1+\mu^{\rm f}(z))}}\bigg)f(x_1-1,x_2)\bigg|^2
\\
\nonumber
&\quad+2\sum_{x\in\Vc^{\rm f}}m^{\rm f}(x)\bigg|\frac{\rmi}{2}\Big(e^{1/2}\langle Q-1/2\rangle^{1+\epsilon}\otimes\un_{\Vc_2^{\rm f}}\Big)
\\
\label{e:part2}
&\hspace*{+0.5cm}\times\sum_{z_2\sim x_2}\delta_{z_1=x_1}\bigg(\frac{\Ec^{\rm f}_2(x_2,z_2)(1+\varepsilon^{\rm f}(x,z))}{ m(x)\sqrt{(1+\mu^{\rm f}(x))(1+\mu^{\rm f}(z))}}
\\
\nonumber
&\hspace*{+0.5cm}-\frac{\Ec^{\rm f}_2(x_2,z_2)(1+\varepsilon^{\rm f}((x_1-1,x_2),z))}{m^{\rm f}_1(x_1-1)m_2\sqrt{(1+\mu^{\rm f}(x_1-1,x_2)(1+\mu^{\rm f}(z))}}\bigg)f(x_1-1,x_2)\bigg|^2.\end{align}
 We bound \eqref{e:part1} as follows:
\begin{align}
\nonumber
&\eqref{e:part1}
\nonumber
\leq4\sum_{x\in\Vc^{\rm f}}m^{\rm f}(x)\bigg|\frac{\rmi}{2}\Big(e \langle Q-1/2\rangle^{1+\epsilon}\otimes\un_{\Vc_2^{\rm f}}\Big)\frac{1}{m_2}
\\
&\label{e:partA}\times\bigg(\frac{\sqrt{1+\mu^{\rm f}(x_1-1,x_2)}(1+\varepsilon^{\rm f}(x,(x_1+1,x_2))}{ \sqrt{(1+\mu^{\rm f}(x))(1+\mu^{\rm f}(x_1+1,x_2))(1+\mu^{\rm f}(x_1-1,x_2))}}
\\
\nonumber
&\hspace*{+0.5cm}-\frac{\sqrt{1+\mu^{\rm f}(x_1+1,x_2)}(1+\varepsilon^{\rm f}((x_1-1,x_2),(x_1,x_2)))}{\sqrt{(1+\mu^{\rm f}(x_1-1,x_2)(1+\mu^{\rm f}(x))(1+\mu^{\rm f}(x_1+1,x_2))}}\bigg)f(x_1-1,x_2)\bigg|^2
\\
\nonumber
&\quad+4\sum_{x\in\Vc^{\rm f}}m^{\rm f}(x)\bigg|\frac{\rmi}{2}\Big( \langle Q-1/2\rangle^{1+\epsilon}\otimes\un_{\Vc_2^{\rm f}}\Big)\frac{1}{m_2}
\\
&\label{e:partB}\hspace*{+0.5cm}\times\bigg(\frac{\sqrt{1+\mu^{\rm f}(x_1-2,x_2)}(1+\varepsilon^{\rm f}(x,(x_1-1,x_2)))}{ \sqrt{(1+\mu^{\rm f}(x))(1+\mu^{\rm f}(x_1-1,x_2))(1+\mu^{\rm f}(x_1-2,x_2))}}
\\
\nonumber
&\hspace*{+0.5cm}-\frac{\sqrt{1+\mu^{\rm f}(x_1-1,x_2)}(1+\varepsilon^{\rm f}((x_1-1,x_2),(x_1-2,x_2)))}{\sqrt{(1+\mu^{\rm f}(x_1-1,x_2)(1+\mu^{\rm f}(x_1-2,x_2))(1+\mu^{\rm f}(x_1-1,x_2))}}\bigg)f(x_1-1,x_2)\bigg|^2,\end{align}
Now, we concentrate on \eqref{e:partA} and in the same way, we deal with \eqref{e:partB}. Since the assertions $(H2)$ and $(H3)$ hold true then there exists an integer $c,$ such that
\begin{align*}&4\sum_{x\in\Vc^{\rm f}}m^{\rm f}(x)\bigg|\frac{\rmi}{2}\Big(e \langle Q-1/2\rangle^{1+\epsilon}\otimes\un_{\Vc_2^{\rm f}}\Big)
\\
&\quad\times\frac{1}{m_2}\bigg(\frac{\sqrt{1+\mu^{\rm f}(x_1-1,x_2)}(1+\varepsilon^{\rm f}(x,(x_1+1,x_2))}{ \sqrt{(1+\mu^{\rm f}(x))(1+\mu^{\rm f}(x_1+1,x_2))(1+\mu^{\rm f}(x_1-1,x_2))}}
\\
&\quad-\frac{\sqrt{1+\mu^{\rm f}(x_1+1,x_2)}(1+\varepsilon^{\rm f}((x_1-1,x_2),(x_1,x_2)))}{\sqrt{(1+\mu^{\rm f}(x_1-1,x_2)(1+\mu^{\rm f}(x))(1+\mu^{\rm f}(x_1+1,x_2))}}\bigg)f(x_1-1,x_2)\bigg|^2\end{align*}
\begin{align*}
&=4\sum_{x\in\Vc^{\rm f}}m^{\rm f}(x)\bigg|\frac{\rmi}{2}\Big(e \langle Q-1/2\rangle^{1+\epsilon}\otimes\un_{\Vc_2^{\rm f}}\Big)\frac{1}{m_2}
\\
&\quad\times\bigg(\frac{(\mu^{\rm f}(x_1-1,x_2)-\mu^{\rm f}(x_1+1,x_2))\varepsilon^{\rm f}(x,(x_1+1,x_2))}{ (\sqrt{1+\mu^{\rm f}(x_1-1,x_2)}+\sqrt{1+\mu^{\rm f}(x_1+1,x_2)})}
\\
&\quad\times\frac{1}{\sqrt{(1+\mu^{\rm f}(x))(1+\mu^{\rm f}(x_1+1,x_2))(1+\mu^{\rm f}(x_1-1,x_2))}}
\\
&\quad+\frac{\sqrt{1+\mu^{\rm f}(x_1+1,x_2)}(\varepsilon^{\rm f}(x,(x_1+1,x_2))-\varepsilon^{\rm f}((x_1-1,x_2),x)}{ \sqrt{(1+\mu^{\rm f}(x))(1+\mu^{\rm f}(x_1+1,x_2))(1+\mu^{\rm f}(x_1-1,x_2))}}
\\
&\quad+\frac{\mu^{\rm f}(x_1-1,x_2)-\mu^{\rm f}(x_1+1,x_2)}{ (\sqrt{1+\mu^{\rm f}(x_1-1,x_2)}+\sqrt{1+\mu^{\rm f}(x_1+1,x_2)})}
\\
&\quad\times\frac{1}{\sqrt{(1+\mu^{\rm f}(x))(1+\mu^{\rm f}(x_1+1,x_2))(1+\mu^{\rm f}(x_1-1,x_2))}}\bigg)f(x_1-1,x_2)\bigg|^2
\\
&\leq  c\|f\|_{\ell^2(\Vc^{\rm f},m^{\rm f})}^2.\end{align*}
In the same way, we treat \eqref{e:part2} and $\|\langle\Lambda^{\rm f}\rangle^{\epsilon}[W^{\rm f}(\cdot),\Ac_{\Gc^{\rm f}}]f\|^2_{\ell^2(\Vc^{\rm f},m^{\rm f})}$.
 By density, there exists $c>0$ such that $\|\langle\Lambda^{\rm f}\rangle^{\epsilon}[\widetilde{\Delta}_{\Gc_{\varepsilon, \mu}^{\rm f}},\Ac_{\Gc^{\rm f}}]f\|^2_{\ell^2(\Vc^{\rm f},m^{\rm f})}\leq c\|f\|^2_{\ell^2(\Vc^{\rm f},m^{\rm f})}$. Finally,  by applying  \cite[Proposition 7.5.7]{ABG} where the hypotheses are verified in Proposition \ref{p:hypothesis}, we find the result. \qed

%
%
%
%
%
%
%
%
%
%
%
%


\subsection{The cusp side: Radial metric perturbation}    We recall that
\[\widetilde{\Delta}_{\Gc_{\varepsilon, \mu}^{\rm c}}f(x):=T_{m_{\mu}\rightarrow m}\Delta_{\Gc_{\varepsilon, \mu}^{\rm c}}T^{-1}_{m_{\mu}\rightarrow m}f(x),\quad \rm{ for\, all\, }\, f\in\Cc_c(\Vc^{\rm c}). \] 
We first deal with the question of the essential spectrum.

\begin{proposition}\label{com.cusp}
Let $V^{\rm c}:\Vc^{\rm c}\to \R$ be a function obeying $V^{\rm c}(x)\to0$ if $|x|\to\infty$. We assume that \eqref{e:tend0} holds true 
then $\widetilde{\Delta}_{\Gc_{\varepsilon, \mu}^{\rm c}}-\Delta_{\Gc^{\rm c}}\in \Kc(\ell^2(\Vc^{\rm c},m^{\rm c}))$. In particular,
\begin{enumerate}
\item[(1)]$\Dc(\Delta_{\Gc_{\varepsilon, \mu}^{\rm c}}+V(\cdot))=\Dc(T_{m_{\mu}\rightarrow m}^{-1}\Delta_{\Gc}T_{m_{\mu}\rightarrow m})$,
\item[(2)]$\Delta_{\Gc_{\varepsilon, \mu}^{\rm c}}+V(\cdot)$ is essentially self-adjoint on $\Cc_c(\Vc)$,
\item[(3)]$\sigma_{\rm ess}(\Delta_{\Gc_{\varepsilon, \mu}^{\rm c}}+V(\cdot))=\sigma_{\rm ess}(\Delta_{\Gc})$.
\end{enumerate}
\end{proposition}
\proof Let $f\in\Cc_c(\Vc^{\rm c})$,
we have
\begin{align*}&|\langle f, (\widetilde{\Delta}_{\Gc_{\varepsilon, \mu}^{\rm c}}-\Delta_{\Gc^{\rm c}}) f\rangle_{\ell^2(\Vc^{\rm c},m^{\rm c})}|=|\sum_{x\in\Vc^{\rm c}} m^{\rm c}(x)(\widetilde{\Delta}_{\Gc_{\varepsilon, \mu}^{\rm c}}-\Delta_{\Gc^{\rm c}}) f(x)\overline{f(x)}|
\\
&\leq\sum_{x\in\Vc^{\rm c}} m^{\rm c}(x)\frac{1}{m^{\rm c}(x)}\sum_{z_1\sim x_1}\frac{\varepsilon^{\rm c}(x,(z_1,x_2))}{\sqrt{(1+\mu^{\rm c}(x))(1+\mu^{\rm c}(z_1,x_2))}}
\\
&\hspace*{+1cm}\times\Ec^{\rm c}_1(x_1,z_1)|f(x)|^2
\\
&\quad+1/2\sum_{x\in\Vc^{\rm c}} m^{\rm c}(x)\frac{1}{m^{\rm c}(x)}\sum_{z_1\sim x_1}\frac{\varepsilon^{\rm c}(x,(z_1,x_2))}{\sqrt{(1+\mu^{\rm c}(x))(1+\mu^{\rm c}(z_1,x_2))}}
\\
&\hspace*{+1cm}\times\Ec^{\rm c}(x,z)(|f(z)|^2+|f(x)|^2)
\\
&\quad+\sum_{x\in\Vc^{\rm c}} m^{\rm c}(x)\frac{1}{m^{\rm c}(x)}\sum_{z_1\sim x_1}\left|\frac{\mu^{\rm c}(x)+\mu^{\rm c}(z_1,x_2)+\mu^{\rm c}(x)\mu^{\rm c}(z_1,x_2)}{\sqrt{(1+\mu^{\rm c}(x))(1+\mu^{\rm c}(z_1,x_2))}}\right|
\\
&\hspace*{+1cm}\times\frac{1}{\sqrt{1+\mu^{\rm c}(x)}+\sqrt{1+\mu^{\rm c}(z_1,x_2)}}\Ec^{\rm c}_1(x_1,z_1)|f(x)|^2
\\
&\quad+1/2 \sum_{x\in\Vc^{\rm c}} m^{\rm c}(x)\frac{1}{m^{\rm c}(x)}\sum_{z_1\sim x_1}\left|\frac{\mu^{\rm c}(x)+\mu^{\rm c}(z_1,x_2)+\mu^{\rm c}(x)\mu^{\rm c}(z_1,x_2)}{\sqrt{(1+\mu^{\rm c}(x))(1+\mu^{\rm c}(z_1,x_2))}}\right|
\\
&\hspace*{+1cm}\times\frac{1}{\sqrt{1+\mu^{\rm c}(x)}+\sqrt{1+\mu^{\rm c}(z_1,x_2)}}\Ec^{\rm c}_1(x_1,z_1)|f(z_1,x_2)|^2
\\
&\quad+1/2 \sum_{x\in\Vc^{\rm c}} m^{\rm c}(x)\frac{1}{m^{\rm c}(x)}\sum_{z_1\sim x_1}\left|\frac{\mu^{\rm c}(x)+\mu^{\rm c}(z_1,x_2)+\mu^{\rm c}(x)\mu^{\rm c}(z_1,x_2)}{\sqrt{(1+\mu^{\rm c}(x))(1+\mu^{\rm c}(z_1,x_2))}}\right|
\\
&\hspace*{+1cm}\times\frac{1}{\sqrt{1+\mu^{\rm c}(x)}+\sqrt{1+\mu^{\rm c}(z_1,x_2)}}\Ec^{\rm c}_1(x_1,z_1)|f(x)|^2
\\
&\leq 2\langle f, (\deg_3(\cdot)+\deg_4(\cdot)+|W^{\rm c}(\cdot)|) f\rangle, \end{align*}
 where
 \[\deg_3(x):=\frac{1}{m^{\rm c}(x)}\sum_{z_1\in\Vc^{\rm c}_1}\frac{\varepsilon^{\rm c}(x,(z_1,x_2))}{\sqrt{(1+\mu^{\rm c}(x))(1+\mu^{\rm c}(z_1,x_2))}}\Ec^{\rm c}_1(x_1,z_1)\]
and
\begin{align*}\deg_4(x)&:=\frac{1}{m^{\rm c}(x)}\sum_{z_1\in\Vc^{\rm c}_1}\left|\frac{\mu^{\rm c}(x)+\mu^{\rm c}(z_1,x_2)+\mu^{\rm c}(x)\mu^{\rm c}(z_1,x_2)}{\sqrt{(1+\mu^{\rm c}(x))(1+\mu^{\rm c}(z_1,x_2))}(\sqrt{1+\mu^{\rm c}(x)}+\sqrt{1+\mu^{\rm c}(z_1,x_2)})}\right|
\\
&\hspace*{+4cm}\times\Ec^{\rm c}_1(x_1,z_1).\end{align*}
We have
 \begin{align*}|\deg_3(x)|&\leq\sup_{z_1\sim x_1}\left|\frac{\varepsilon^{\rm c}(x,(z_1,x_2))}{m_2\sqrt{(1+\mu^{\rm c}(x))(1+\mu^{\rm c}(z_1,x_2))}}\right|\deg_{\Gc_1^{\rm c}}(x).\end{align*}
Since $V_2$ is a finite set and for all $x_2\in\Vc_2$, $\frac{\varepsilon^{\rm c}((x_1,x_2),(z_1,x_2))}{\sqrt{(1+\mu^{\rm c}(x_1,x_2))(1+\mu^{\rm c}(z_1,x_2))}}\to0$  when $x_1,z_1\to\infty$ and since $\deg_{\Gc_1^{\rm c}}(\cdot)$ is bounded then $\deg_3(\cdot)$ is compact. In the same way, using that $\forall x_2\in\Vc_2,~\frac{\mu^{\rm c}(x)+\mu^{\rm c}(z_1,x_2)+\mu^{\rm c}(x)\mu^{\rm c}(z_1,x_2)}{\sqrt{(1+\mu^{\rm c}(x_1,x_2))(1+\mu^{\rm c}(z_1,x_2))}(\sqrt{1+\mu^{\rm c}(x)}+\sqrt{1+\mu^{\rm c}(z_1,x_2)})}\to 0$ if $x_1,z_1\to\infty$, we obtain the compactness of $\deg_4(\cdot)$. 

Now, we will show that $W^{\rm c}(\cdot)\in\Kc(\ell^2(\Vc^{\rm c},m^{\rm c})$. For all  $x\in\Vc^{\rm c},$  we have
\begin{align*}&|W^{\rm c}(x)|=\big|\frac{1}{m^{\rm c}(x)}\sum_{z\sim x}\left(\frac{\mu^{\rm c}(z)-\mu^{\rm c}(x)}{(1+\mu^{\rm c}(x))\sqrt{1+\mu^{\rm c}(z)}(\sqrt{1+\mu^{\rm c}(z)}+\sqrt{1+\mu^{\rm c}(x)})}\right)
\\
&\quad\times(1+\varepsilon^{\rm c}(x,z))\Ec^{\rm c}(x,z)\big|
\\
&\leq\sup_{z_1\sim x_1}(1+\varepsilon^{\rm c}(x,(z_1,x_2)))\deg_{\Gc_1^{\rm c}}(x)
\\
&\quad \times
\left|\frac{\mu^{\rm c}(z_1,x_2)-\mu^{\rm c}(x)}{(1+\mu^{\rm c}(x))\sqrt{1+\mu^{\rm c}(z_1,x_2)}(\sqrt{1+\mu^{\rm c}(z_1,x_2)}+\sqrt{1+\mu^{\rm c}(x)})}\right|\end{align*}
Since $\Vc_2$ is a finite set and $\forall x_2\in\Vc_2$, $\varepsilon^{\rm c}((x_1,x_2),(z_1,x_2))(\mu^{\rm c}(z_1,x_2)-\mu^{\rm c}(x_1,x_2))\to0$ when $x_1,z_1\to\infty$, and since $\deg_{\Gc_1^{\rm c}}(\cdot)$ is bounded and since $V^{\rm c}(\cdot)$ is a compact perturbation. Then, $\widetilde{\Delta}_{\Gc_{\varepsilon, \mu}^{\rm c}}-\Delta_{\Gc^{\rm c}}$ is a compact operator.  The points (1) and (2) follow from Theorem \cite[Theorem XIII.14]{RS} and (3) from the Weyl's Theorem. \qed

In order to go into the Mourre theory, we construct the conjugate operator:
\begin{align}\label{Acusp}
\Ac_{\Gc^{\rm c}}&:=\Ac_{m_1}\otimes P^{\rm le},\end{align} with \begin{align*}\Ac_{m_1}&:=T_{1\rightarrow m_1}\Ac_{\N}T^{-1}_{1\rightarrow m_1}
\\
&:=\frac{\rmi}{2}\left(e^{1/2}(Q-1/2)U-e^{-1/2}(Q+1/2)U^*\right).
\end{align*}
It is self-adjoint and essentially self-adjoint on $\Cc^c(\Vc^c)$ by Lemma \ref{l:ANessaa}. Because of the projection in \eqref{Acusp}, we restrict to  radial perturbations.
\begin{definition}\label{d:rad}
The perturbations $V^{\rm c}$, $\mu$ and $\varepsilon$ are called $\textbf{radial}$ if they do not depend on the second variable, i.e., For all $(x_1,x_2),(z_1,z_2)\in\Vc^{\rm c}$, we have $V^{\rm c}(x_1,x_2)=V^{\rm c}(x_1,z_2)$, $\mu(x_1,x_2)=\mu(x_1,z_2)$
and $\varepsilon((x_1,x_2),(z_1,z_2))=\varepsilon((x_1,x_2),(z_1,x_2))$.
\end{definition}

We turn to Mourre estimate.

\begin{proposition}\label{mourcusp2}Let $\Gc_{\varepsilon, \mu}^{\rm c}$ a graph satisfies a condition \eqref{e:tend0}. Suppose that $V^{\rm c}:\Vc^{\rm c}\to \R$, $\varepsilon$ and $\mu$ are radial and  assume (H1), (H2), and (H3)
 and 
 $V^{\rm c}(x)\to0$ if $|x|\to\infty$. Then $\Delta_{\Gc_{\varepsilon, \mu}^{\rm c}}+V^{\rm c}(\cdot)\in\Cc^{1,1}(\Ac_{\Gc_{\varepsilon, \mu}^{\rm c}})$. Moreover, for all compact interval $I\subset (\frac{\alpha}{m_2},\frac{\beta}{m_2})$ there are $c>0$, a compact operator $K$ such that
 \begin{align}
\nonumber
 E_I(\Delta_{\Gc_{\varepsilon, \mu}^{\rm c}}+V^{\rm c}(\cdot))[\Delta_{\Gc_{\varepsilon, \mu}^{\rm c}}+V^{\rm c}(\cdot), \rmi \Ac_{\Gc_{\varepsilon, \mu}^{\rm c}}]_\circ E_I(\Delta_{\Gc_{\varepsilon, \mu}^{\rm c}}+V^{\rm c}(\cdot))&
 \\
 \label{mourrC}
 & \hspace*{-2cm}\geq c E_I(\Delta_{\Gc_{\varepsilon, \mu}^{\rm c}}+V^{\rm c}(\cdot)) + K,\end{align} in the form sense.
\end{proposition}
\proof The Proposition \ref{C1,1cusp} and Lemma \ref{l:Vcusp} gives that $\Delta_{\Gc_{\varepsilon, \mu}^{\rm c}}\in\Cc^{1,1}(\Ac_{\Gc^{\rm c}})$. Since $\widetilde{\Delta}_{\Gc_{\varepsilon, \mu}^{\rm c}}-\Delta_{\Gc^{\rm c}}$ is a compact operator by Proposition \ref{com.cusp}, thanks to \eqref{e:mourrec} and by \cite[Theorem 7.2.9]{ABG} we obtain \eqref{mourrC}.\qed

 We turn to series of  Lemmata.
To be able to apply the \cite[Proposition 7.5.7]{ABG}, we check the next point.
\begin{proposition}\label{c:hypothesis}
Let $\Lambda^{\rm c}:=(Q+1/2)\otimes\un_{\Vc^{\rm c}},$ then $\Lambda^{\rm c}$ satisfies the following assertions:
\begin{enumerate}
\item[(1)] $e^{\rmi\Lambda^{\rm c}t}\Dc(\Delta_{\Gc_{\varepsilon, \mu}^{\rm c}})\subset \Dc(\Delta_{\Gc_{\varepsilon, \mu}^{\rm c}})$  and there exists a finite constant $\rm c$, such that  \[\|e^{\rmi\Lambda^{\rm c}t}\|_{\Bc(\Dc(\Delta_{\Gc_{\varepsilon, \mu}^{\rm c}}))} \leq \rm c, \quad \hbox{for all}\quad t\in\R.\]
\item[(2)] $\Dc(\Lambda^{\rm c})\subset \Dc(\Ac_{\Gc_{\varepsilon, \mu}^{\rm c}})$.
\item[(3)] $(\Lambda^{\rm c})^{-2}(\Ac_{\Gc_{\varepsilon, \mu}^{\rm c}})^2$ extends to a continuous operator in $\Dc(\Delta_{\Gc_{\varepsilon, \mu}^{\rm c}})$.
\end{enumerate}
\end{proposition}
\proof
With the help of the unitary transformation $T_{m_\mu\to m}$, it is enough to prove the result with $\varepsilon=0$ and $\mu=0$.

(1) We have \[[\Delta_{\Gc^{\rm c}},e^{\rmi\Lambda^{\rm c}t }]=[\Delta_{\Gc_1^{\rm c}},e^{\rmi\Lambda^{\rm c}t }]\otimes\frac{1}{m_2}+[\frac{1}{m_1(\cdot)},e^{\rmi\Lambda^{\rm c}t }]\otimes\Delta_2.\] Since $\frac{1}{m_1(\cdot)}$ and $e^{\rmi\Lambda^{\rm c}t }$ commute and since $[\Delta_{\Gc_1^{\rm c}},e^{\rmi\Lambda^{\rm c}t }]$ is uniformly bounded, then  there exists $c>0$ such that for all $f\in\Cc_c(\Vc^{\rm c})$
\[\|(\Delta_{\Gc^{\rm c}}+\rm i)e^{\rmi\Lambda^{\rm c}t }(\Delta_{\Gc^{\rm c}}+\rm i)^{-1}f\|_{\ell^2(\Vc^{\rm c},m^{\rm c})}\leq c\|f\|_{\ell^2(\Vc^{\rm c},m^{\rm c})}.\]
Hence,  there exists $c>0$ such that for all $f\in\Cc_c(\Vc^{\rm c})$
\[\|(\Delta_{\Gc^{\rm c}}+\rm i)e^{\rmi\Lambda^{\rm c}t } f \|_{\ell^2(\Vc^{\rm c},m^{\rm c})}\leq c\|(\Delta_{\Gc^{\rm c}}+\rm i) f \|_{\ell^2(\Vc^{\rm c},m^{\rm c})}.\]
 Since $\Delta_{\Gc^{\rm c}}$ is essentially self-adjoint on $\Cc_c(\Vc^{\rm c})$ then we find the result.

(2) Let $f\in\Cc_c(\Vc^{\rm c})$, by using the relations of Subsection \ref{dec}, we have
\begin{align*}\|\Ac_{\Gc^{\rm c}}f\|^2_{\ell^2(\Vc^{\rm c},m^{\rm c})}
&\leq\frac{1}{2}\sum_{x\in\Vc^{\rm f}}m^{\rm f}(x)\left|e^{1/2}(Q-1/2)U\otimes P^{\rm le}f(x)\right|^2
\\
&\quad+\left|e^{-1/2}(Q+1/2)U^*\otimes P^{\rm le}f(x)\right|^2
\\
&\leq c \sum_{x\in\Vc^{\rm f}}m^{\rm f}(x)\left|(Q+1/2)\otimes P^{\rm le}f(x)\right|^2
\leq c\|\Lambda^{\rm c}f\|^2_{\ell^2(\Vc^{\rm f},m^{\rm f})}.\end{align*}
Since, $\Lambda^{\rm c}$ is essentially self-adjoint on $\Cc_c(\Vc^{\rm c})$. we find the result.

(3) First for all $f\in\Cc_c(\Vc^{\rm c})$, we have
\begin{align*}
&\|(\Lambda^{\rm c})^{-2}(\Ac_{\Gc^{\rm c}})^2f\|^2_{\ell^2(\Vc^{\rm c},m^{\rm c})}=\sum_{(x_1,x_2)\in\Vc^{\rm c}}m^{\rm c}(x_1,x_2)|\Lambda^{-2}(\Ac_{\Gc^{\rm c}})^2f(x,y)|^2
\\
&=\sum_{(x_1,x_2)\in\Vc^{\rm c}}m^{\rm c}(x_1,x_2)\big|\frac{1}{4}\left(\big((Q+1/2)^{-2}(2Q^2+1/2)\big)\otimes P^{\rm le}\right)f(x_1,x_2)
\\
&\quad-\frac{1}{4}e\left(\big((Q+1/2)^{-2}(Q-1/2)(Q-3/2)\big)\otimes P^{\rm le}\right)f(x_1-2,x_2)
\\
&\quad-\frac{1}{4}e^{-1}\left(\big((Q+1/2)^{-2}(Q+1/2)(Q+3/2)\big)\otimes P^{\rm le}\right)f(x_1+2,x_2)\big|^2.\end{align*}
By density, we get $(\Lambda^{\rm c})^{-2}(\Ac_{\Gc^{\rm c}})^2$ is a bounded operator.
Since $\Lambda^{\rm c}$ is a radial operator and $\Delta_{\Gc_1^{\rm c}}$ is bounded then there exists $C>0$ such that, for all $f\in\Cc_c(\Vc^{\rm c})$,
\begin{align*}
\|[\Delta_{\Gc^{\rm c}},(\Lambda^{\rm c})^{-2}(\Ac_{\Gc^{\rm c}})^2]f\|_{\ell^2(\Vc^{\rm c},m^{\rm c})}&
\\
&\hspace*{-4cm}=\|\Big([\Delta_{\Gc_1^{\rm c}}\otimes\frac{1}{m_2},(\Lambda^{\rm c})^{-2}(\Ac_{\Gc^{\rm c}})^2]+[\frac{1}{m_1(\cdot)}\otimes\Delta_{\Gc_2^{\rm c}},(\Lambda^{\rm c})^{-2}(\Ac_{\Gc^{\rm c}})^2]\Big)f\|_{\ell^2(\Vc^{\rm c},m^{\rm c})}
\\
&\hspace*{-4cm}=\|\Big([\Delta_{\Gc_1^{\rm c}}\otimes\frac{1}{m_2},(\Lambda^{\rm c})^{-2}(\Ac_{\Gc^{\rm c}})^2]+(\Lambda^{\rm c})^{-2}[\frac{1}{m_1(\cdot)}\otimes\Delta_{\Gc_2^{\rm c}},(\Ac_{\Gc^{\rm c}})^2]\Big)f\|_{\ell^2(\Vc^{\rm c},m^{\rm c})}
\\
&\hspace*{-4cm}=\|[\Delta_{\Gc_1^{\rm c}}\otimes\frac{1}{m_2},(\Lambda^{\rm c})^{-2}(\Ac_{\Gc^{\rm c}})^2]f\|_{\ell^2(\Vc^{\rm c},m^{\rm c})}\leq C\|f\|_{\ell^2(\Vc^{\rm c},m^{\rm c})}.
\end{align*} We have used  $[\frac{1}{m_1(\cdot)}\otimes\Delta_{\Gc_2^{\rm c}},(\Ac_{\Gc^{\rm c}})^2]=0$ by construction. Conclude by density. \qed

The proof of Proposition \ref{C1,1cusp} is long and technical. For the sake of the reader, we have separated the treatment of the potential $V^{\rm c}$ to present the technical steps.

\begin{lemma}\label{l:Vcusp}Let $V^{\rm c}:\Vc^{\rm c}\to\R$ be a radial function and $(H1)$ holds true, then $[V^{\rm c}(\cdot),\Ac_{\Gc_{\varepsilon, \mu}^{\rm c}}]_\circ\in\Cc^{0,1}(\Ac_{\Gc_{\varepsilon, \mu}^{\rm c}})$. In particular, $V^{\rm c}(\cdot)\in\Cc^{1,1}(\Ac_{\Gc_{\varepsilon, \mu}^{\rm c}})$.
\end{lemma}
\proof Since $V^{\rm c}$ is radial, by a slight abuse of notation,
we have $V^{\rm c}:=V^{\rm c}\otimes\un_{\Vc_2}$. 
We compute the commutator on $\Cc^c(\Vc^{\rm c})$ and get
\begin{align*}
\left[V^{\rm c}(\cdot),\rmi \Ac_{\Gc_{\varepsilon, \mu}^{\rm c}}\right]&=\left(\frac{e^{-1/2}}{2} \left(Q+\frac{1}{2}\right)\left[V^{\rm c},U^*\right] +\frac{e^{1/2}}{2} \left(\frac{1}{2}-Q\right)\left[V^{\rm c},U\right]\right)\otimes P^{\rm le}.\end{align*}
By density, we infer that $\left[V^{\rm c}(\cdot),\rmi \Ac_{\Gc_{\varepsilon, \mu}^{\rm c}}\right]$ extends to a bounded operator and that $V^{\rm c}(\cdot) \in \Cc^1(\Ac_{\Gc_{\varepsilon, \mu}^{\rm c}})$.
Next, there exists $C>0$ so that, for all $f\in\Cc_c(\Vc^{\rm c})$,
 \begin{align*}
 \big\|\langle \Lambda\rangle^{\epsilon}\otimes[V^{\rm c}(\cdot),\rmi \Ac_{\Gc_{\varepsilon, \mu}^{\rm c}}]_\circ f\big\|
&\leq\frac{e^{-1/2}}{2}\big\|\big(\langle Q+1/2\rangle^{\epsilon}(Q+1/2)\left[V^{\rm c},U^*\right]\big)\otimes P^{\rm le}f \big\|
\\
&\quad +\frac{e^{1/2}}{2}\big\|\big(\langle Q+1/2\rangle^{\epsilon} (Q+1/2)\left[V^{\rm c},U\right]\big)\otimes P^{\rm le}f \big\|
\\
&\leq\frac{e^{-1/2}}{2}\big\|\big(\langle \Lambda\rangle^{\epsilon+1}\left[V^{\rm c},U^*\right]\big)\otimes P^{\rm le}f \big\|
\\
&\quad+\frac{e^{1/2}}{2}\big\|\big(\langle \Lambda\rangle^{\epsilon+1}\left[V^{\rm c},U\right]\big)\otimes P^{\rm le}f \big\|
\leq C\|f\|, \quad\rm{ by (H1).}
\end{align*}
Finally, the result follow by applying \cite[Proposition 7.5.7]{ABG} where the hypotheses are verified in Proposition \ref{c:hypothesis}.\qed

Here is the most technical part:
\begin{proposition}\label{C1,1cusp}
Assuming $(H2)$ and $(H3)$, 
 we have ${\Delta}_{\Gc_{\varepsilon, \mu}^{\rm c}}\in \Cc^1(\Ac_{\Gc_{\varepsilon, \mu}^{\rm c}})$. Moreover
$[{\Delta}_{\Gc_{\varepsilon, \mu}^{\rm c}},\Ac_{\Gc_{\varepsilon, \mu}^{\rm c}}]_\circ\in\Cc^{0,1}(\Ac_{\Gc_{\varepsilon, \mu}^{\rm c}})$. In particular, ${\Delta}_{\Gc_{\varepsilon, \mu}^{\rm c}}\in\Cc^{1,1}(\Ac_{\Gc_{\varepsilon, \mu}^{\rm c}})$.
\end{proposition}

\proof We work in $\ell^2(\Vc^{\rm c}, m^{c})$. We first prove that $\widetilde{\Delta}_{\Gc_{\varepsilon, \mu}^{\rm c}}\in\Cc^{1}(\Ac_{\Gc^{\rm c}})$. By the computation below (with $\epsilon=0$), we obtain that there is $c>0$ such that
\[\|[\widetilde{\Delta}_{\Gc_{\varepsilon, \mu}^{\rm c}},\Ac_{\Gc^{\rm c}}] f\|_{\ell^2(\Vc^{\rm c},m^{\rm c})}\leq c \|f\|_{\ell^2(\Vc^{\rm c},m^{\rm c})}, \quad \forall f\in \Cc_c(\Vc^{\rm c}).\]
Using Lemma \ref{l:gec} and \cite[Theorem 6.3.4]{ABG}, this implies that $\widetilde{\Delta}_{\Gc_{\varepsilon, \mu}^{\rm c}}\in\Cc^{1}(\Ac_{\Gc^{\rm c}})$.

We turn to the $\Cc^{0,1}$ property. We assume that $(H2)$ and $(H3)$ are true then
\begin{align}\nonumber
&\|\langle\Lambda^{\rm c}\rangle^{\epsilon}[\widetilde{\Delta}_{\Gc_{\varepsilon, \mu}^{\rm c}},\Ac_{\Gc^{\rm c}}]_\circ f\|_{\ell^2(\Vc^{\rm c},m^{\rm c})}
\leq\sum_{x\in\Vc^{\rm c}}m^{\rm c}(x)|\frac{\rmi}{2}(e^{1/2}\langle Q-1/2\rangle^{1+\epsilon}\otimes  P^{\rm le}
\\
\nonumber
&\hspace*{+0.2cm}\times\big(\frac{1}{m^{\rm c}(x)}\sum_{z\sim x}\Ec^{\rm c}(x,z)\frac{1+\varepsilon^{\rm c}(x,z)}{\sqrt{(1+\mu^{\rm c}(x))(1+\mu^{\rm c}(z))}}
\\
&\label{e:comcusp}\hspace*{+0.2cm}-\frac{1}{m^{\rm c}(x_1-1,x_2)}\sum_{z\sim x}\Ec^{\rm c}((x_1-1,x_2),z)\frac{1+\varepsilon^{\rm c}((x_1-1,x_2),z)}{\sqrt{(1+\mu^{\rm c}(x_1-1,x_2)(1+\mu^{\rm c}(z))}}\big)
\\
\nonumber
&\hspace*{+0.2cm}\times f(x_1-1,x_2)|^2
\\
\nonumber
&\hspace*{+0.1cm}+\sum_{x\in\Vc^{\rm c}}m^{\rm c}(x)|\frac{\rmi}{2}(e^{1/2}\langle Q-1/2\rangle^{1+\epsilon}\otimes  P^{\rm le})
\\
\nonumber
&\hspace*{+0.2cm}\times\big(\frac{1}{m^{\rm c}(x_1-1,x_2)}\sum_{z\sim x}\Ec^{\rm c}((x_1-1,x_2),z)\frac{1+\varepsilon^{\rm c}((x_1-1,x_2),z)}{\sqrt{(1+\mu^{\rm c}(x_1-1,x_2)(1+\mu^{\rm c}(z))}}f(z)
\\
\nonumber
&\hspace*{+0.2cm}-\frac{1}{m^{\rm c}(x)}\sum_{z\sim x}\Ec^{\rm c}(x,z)\frac{1+\varepsilon^{\rm c}(x,z)}{\sqrt{(1+\mu^{\rm c}(x))(1+\mu^{\rm c}(z))}}f(z_1-1,z_2)\big)|^2
\\
\nonumber
&\hspace*{+0.1cm}+\sum_{x\in\Vc^{\rm c}}m^{\rm c}(x)|\frac{\rmi}{2}(e^{-1/2}\langle Q+1/2\rangle^{1+\epsilon}\otimes  P^{\rm le})
\\
\nonumber
&\hspace*{+0.2cm}\times\big(\frac{1}{m^{\rm c}(x_1+1,x_2)}\sum_{z\sim x}\Ec^{\rm c}((x_1+1,x_2),z)\frac{1+\varepsilon^{\rm c}((x_1+1,x_2),z)}{\sqrt{(1+\mu^{\rm c}(x_1+1,x_2)(1+\mu^{\rm c}(z))}}
\\
\nonumber
&\hspace*{+0.2cm}-\frac{1}{m^{\rm c}(x)}\sum_{z\sim x}\Ec^{\rm c}(x,z)\frac{1+\varepsilon^{\rm c}(x,z)}{\sqrt{(1+\mu^{\rm c}(x))(1+\mu^{\rm c}(z))}}\big)f(x_1+1,x_2)|^2
\\
\nonumber
&\hspace*{+0.1cm}+\sum_{x\in\Vc^{\rm c}}m^{\rm c}(x)|\frac{\rmi}{2}(e^{-1/2}\langle Q+1/2\rangle^{1+\epsilon}\otimes  P^{\rm le})
\\
\nonumber
&\hspace*{+0.2cm}\times\big(\frac{1}{m^{\rm c}(x)}\sum_{z\sim x}\Ec^{\rm c}(x,z)\frac{1+\varepsilon^{\rm c}(x,z)}{\sqrt{(1+\mu^{\rm c}(x))(1+\mu^{\rm c}(z))}}f(z_1+1,z_2)
\\
\nonumber
&\hspace*{+0.2cm}-\frac{1}{m^{\rm c}(x_1+1,x_2)}\sum_{z\sim x}\Ec^{\rm c}((x_1+1,x_2),z)\frac{1+\varepsilon^{\rm c}((x_1+1,x_2),z)}{\sqrt{(1+\mu^{\rm c}(x_1+1,x_2)(1+\mu^{\rm c}(z))}}f(z_1,z_2)\big)|^2
\\
\nonumber&+\|\langle\Lambda^{\rm c}\rangle^{\epsilon}[W^{\rm c},\Ac_{\Gc^{\rm c}}]_\circ f\|_{\ell^2(\Vc^{\rm c},m^{\rm c})}.\end{align}
 We treat the first term of $\|\langle\Lambda^{\rm c}\rangle^{\epsilon}[\widetilde{\Delta}_{\Gc_{\varepsilon, \mu}^{\rm c}},\Ac_{\Gc^{\rm c}}]_\circ f\|_{\ell^2(\Vc^{\rm c},m^{\rm c})}$ in \eqref{e:comcusp}
\begin{align}\nonumber&\sum_{x\in\Vc^{\rm c}}m^{\rm c}(x)\big|\frac{\rmi}{2}\left(e^{1/2}\langle Q-1/2\rangle^{1+\epsilon}\otimes  P^{\rm le}\right) 
\\
\nonumber
&\hspace*{+1cm}\times\sum_{z\sim x}\big(\frac{\Ec^{\rm c}(x,z)(1+\varepsilon^{\rm c}(x,z))}{ m^{\rm c}(x)\sqrt{(1+\mu^{\rm c}(x))(1+\mu^{\rm c}(z))}}
\\
\nonumber
&\hspace*{+1cm}-\frac{\Ec^{\rm c}((x_1-1,x_2),z)(1+\varepsilon^{\rm c}((x_1-1,x_2),z))}{m^{\rm c}(x_1-1,x_2)\sqrt{(1+\mu^{\rm c}(x_1-1,x_2)(1+\mu^{\rm c}(z))}}\big)f(x_1-1,x_2)\big|^2
\\
\nonumber
&\leq2\sum_{x\in\Vc^{\rm c}}m^{\rm c}(x)\big|\frac{\rmi}{2}\left(e^{1/2}\langle Q-1/2\rangle^{1+\epsilon}\otimes  P^{\rm le}\right) 
\\
&\label{e:Part1}\hspace*{+1cm}\times\sum_{z_1\sim x_1}\delta_{z_2=x_2}\big(\frac{\Ec^{\rm c}_1(x_1,z_1)(1+\varepsilon^{\rm c}(x,z))}{ m(x)\sqrt{(1+\mu^{\rm c}(x))(1+\mu^{\rm c}(z))}}
\\
\nonumber
&\hspace*{+1cm}-\frac{\Ec^{\rm c}_1(x_1-1,z_1)(1+\varepsilon^{\rm c}((x_1-1,x_2),z))}{m^{\rm c}(x_1-1,x_2)\sqrt{(1+\mu^{\rm c}(x_1-1,x_2)(1+\mu^{\rm c}(z))}}\big)f(x_1-1,x_2)|^2\end{align}

\begin{align}\nonumber
&+2\sum_{x\in\Vc^{\rm c}}m^{\rm c}(x)\big|\frac{\rmi}{2}\left(e^{1/2}\langle Q-1/2\rangle^{1+\epsilon}\otimes  P^{\rm le}\right) 
\\
&\label{e:Part2}\hspace*{+1cm}\times\sum_{z_2\sim x_2}\delta_{z_1=x_1}\big(\frac{\Ec_2(x_2,z_2)(1+\varepsilon^{\rm c}(x,z))}{ m(x)\sqrt{(1+\mu^{\rm c}(x))(1+\mu^{\rm c}(z))}}
\\
\nonumber
&\hspace*{+1cm}-\frac{\Ec_2(x_2,z_2)(1+\varepsilon^{\rm c}((x_1-1,x_2),z))}{m^{\rm c}_1(x_1-1)m_2\sqrt{(1+\mu^{\rm c}(x_1-1,x_2)(1+\mu^{\rm c}(z))}}\big)f(x_1-1,x_2)|^2.\end{align}
We focus on \eqref{e:Part1}.
\begin{align}\eqref{e:Part1}\nonumber&\leq4\sum_{x\in\Vc^{\rm c}}m^{\rm c}(x)\big|\frac{\rmi}{2}\left(e \langle Q-1/2\rangle^{1+\epsilon}\otimes  P^{\rm le}\right) 
\\
&\label{e:PartA}\hspace*{+0.5cm}\times\frac{1}{m_2}\big(\frac{\sqrt{1+\mu^{\rm c}(x_1-1,x_2)}(1+\varepsilon^{\rm c}(x,(x_1+1,x_2))}{ \sqrt{(1+\mu^{\rm c}(x))(1+\mu^{\rm c}(x_1+1,x_2))(1+\mu^{\rm c}(x_1-1,x_2))}}
\\
\nonumber
&\hspace*{+1cm}-\frac{\sqrt{1+\mu^{\rm c}(x_1+1,x_2)}(1+\varepsilon^{\rm c}((x_1-1,x_2),(x_1,x_2)))}{\sqrt{(1+\mu^{\rm c}(x_1-1,x_2)(1+\mu^{\rm c}(x))(1+\mu^{\rm c}(x_1+1,x_2))}}\big)f(x_1-1,x_2)\big|^2
\\
\nonumber
&+4\sum_{x\in\Vc^{\rm c}}m^{\rm c}(x)\big|\frac{\rmi}{2}\left( \langle Q-1/2\rangle^{1+\epsilon}\otimes  P^{\rm le}\right) 
\\
\nonumber
&\hspace*{+0.5cm}\times\frac{1}{m_2}\big(\frac{\sqrt{1+\mu^{\rm c}(x_1-2,x_2)}(1+\varepsilon^{\rm c}(x,(x_1-1,x_2)))}{ \sqrt{(1+\mu^{\rm c}(x))(1+\mu^{\rm c}(x_1-1,x_2))(1+\mu^{\rm c}(x_1-2,x_2))}}
\\
&\label{e:PartB}\hspace*{+1cm}-\frac{\sqrt{1+\mu^{\rm c}(x_1-1,x_2)}(1+\varepsilon^{\rm c}((x_1-1,x_2),(x_1-2,x_2)))}{\sqrt{(1+\mu^{\rm c}(x_1-1,x_2)(1+\mu^{\rm c}(x_1-2,x_2))(1+\mu^{\rm c}(x_1-1,x_2))}}\big)
\\
\nonumber
&\hspace*{+1cm}\times f(x_1-1,x_2)\big|^2.\end{align}
Now, we concentrate on \eqref{e:PartA}. \eqref{e:PartB} can be done in the same way. Since the assertions $(H2)$ and $(H3)$ hold true then there exists an integer $c,$ such that \begin{align*}&4\sum_{x\in\Vc^{\rm c}}m^{\rm c}(x)\big|\frac{\rmi}{2}\left(e (Q-1/2)^{1+\epsilon}\otimes  P^{\rm le}\right) 
\\
&\hspace*{+0.5cm}\times\frac{1}{m_2}\big(\frac{\sqrt{1+\mu^{\rm c}(x_1-1,x_2)}(1+\varepsilon^{\rm c}(x,(x_1+1,x_2))}{ \sqrt{(1+\mu^{\rm c}(x))(1+\mu^{\rm c}(x_1+1,x_2))(1+\mu^{\rm c}(x_1-1,x_2))}}
\\
&\hspace*{+1cm}-\frac{\sqrt{1+\mu^{\rm c}(x_1+1,x_2)}(1+\varepsilon^{\rm c}((x_1-1,x_2),(x_1,x_2)))}{\sqrt{(1+\mu^{\rm c}(x_1-1,x_2)(1+\mu^{\rm c}(x))(1+\mu^{\rm c}(x_1+1,x_2))}}\big)f(x_1-1,x_2)\big|^2
\\
&=4\sum_{x\in\Vc^{\rm c}}m^{\rm c}(x)\big|\frac{\rmi}{2}\left(e \langle Q-1/2\rangle^{1+\epsilon}\otimes  P^{\rm le}\right)\frac{1}{m_2} 
\\
&\hspace*{+0.5cm}\times\big(\frac{(\mu^{\rm c}(x_1-1,x_2)-\mu^{\rm c}(x_1+1,x_2))\varepsilon^{\rm c}(x,(x_1+1,x_2))}{ (\sqrt{1+\mu^{\rm c}(x_1-1,x_2)}+\sqrt{1+\mu^{\rm c}(x_1+1,x_2)})}
\\
&\hspace*{+0.7cm} \times\frac{1}{\sqrt{(1+\mu^{\rm c}(x))(1+\mu^{\rm c}(x_1+1,x_2))(1+\mu^{\rm c}(x_1-1,x_2))}}\big)\end{align*}
\begin{align*}&+4\sum_{x\in\Vc^{\rm c}}m^{\rm c}(x)\big|\frac{\rmi}{2}\left(e \langle Q-1/2\rangle^{1+\epsilon}\otimes  P^{\rm le}\right)\frac{1}{m_2}
\\
&\hspace*{+1cm}\times\big(\frac{\sqrt{1+\mu^{\rm c}(x_1+1,x_2)}(\varepsilon^{\rm c}(x,(x_1+1,x_2))-\varepsilon^{\rm c}((x_1-1,x_2),x)}{ \sqrt{(1+\mu^{\rm c}(x))(1+\mu^{\rm c}(x_1+1,x_2))(1+\mu^{\rm c}(x_1-1,x_2))}}
\\
&\hspace*{+1cm}+\frac{\mu^{\rm c}(x_1-1,x_2)-\mu^{\rm c}(x_1+1,x_2)}{ (\sqrt{1+\mu^{\rm c}(x_1-1,x_2)}+\sqrt{1+\mu^{\rm c}(x_1+1,x_2)})}
\\
&\hspace*{+0.7cm}\times\frac{1}{\sqrt{(1+\mu^{\rm c}(x))(1+\mu^{\rm c}(x_1+1,x_2))(1+\mu^{\rm c}(x_1-1,x_2))}}\big)f(x_1-1,x_2)\big|^2
\\
&\leq  c\|f\|_{\ell^2(\Vc^{\rm c},m^{\rm c})}^2,\end{align*}
and in the same way, we deal with \eqref{e:Part2} and $\|\langle\Lambda^{\rm c}\rangle^{\epsilon}[W^{\rm c},\Ac_{\Gc^{\rm c}}]_\circ f\|_{\ell^2(\Vc^{\rm c},m^{\rm c})}$.  By density, we have proven that there exists $c>0$ such that
\[\|\langle\Lambda^{\rm c}\rangle^{\epsilon}[\widetilde{\Delta}_{\Gc_{\varepsilon, \mu}^{\rm c}},\Ac_{\Gc^{\rm c}}]_\circ f\|^2_{\ell^2(\Vc^{\rm c},m^{\rm c})}\leq c\|f\|^2_{\ell^2(\Vc^{\rm c},m^{\rm c})}.\]
 Finally, by applying \cite[Proposition 7.5.7]{ABG} with $\Gc:= \Dc(\Delta_{\Gc^{\rm c}_{\varepsilon, \mu}})$ where the hypotheses are verified in Proposition \ref{p:hypothesis}, we find the result. \qed


\end{document}